\documentclass[a4paper]{article}
\usepackage[dvips]{graphicx}
\usepackage{amsmath}
\usepackage{amsfonts}
\usepackage{amssymb}
\newtheorem{theorem}{Theorem}[subsection]

\newtheorem{corollary}[theorem]{Corollary}

\newtheorem{definition}[theorem]{Definition}
\newtheorem{example}[theorem]{Example}

\newtheorem{lemma}[theorem]{Lemma}
\newtheorem{notation}[theorem]{Notation}

\newtheorem{proposition}[theorem]{Proposition}
\newtheorem{remark}[theorem]{Remark}

\newenvironment{proof}[1][Proof]{\textbf{#1.} }{\ \rule{0.5em}{0.5em}}

\begin{document}

\title{Schensted-type correspondences and plactic monoids for types $B_{n}$ and
$D_{n}$}
\author{Cedric Lecouvey\\lecouvey@math.unicaen.fr}
\date{\ }
\maketitle
\begin{abstract}
We use Kashiwara's theory of crystal bases to study plactic monoids for
$U_{q}(so_{2n+1})$ and $U_{q}(so_{2n})$. Simultaneously we describe a
Schensted type correspondence in the crystal graphs of tensor powers of vector
and spin representations and we derive a Jeu de Taquin for type $B$ from the
Sheats sliding algorithm.
\end{abstract}

\section{ Introduction}

The Schensted correspondence based on the bumping algorithm yields a bijection
between words $w$ of length $l$ on the ordered alphabet $\mathcal{A}%
_{n}=\{1\prec2\prec\cdot\cdot\cdot\prec n\}$ and pairs $(P^{A}(w),Q^{A}(w))$
of tableaux of the same shape containing $l$ boxes where $P^{A}(w)$ is a
\ semi-standard Young tableau on $\mathcal{A}_{n}$ and $Q^{A}(w)$ is a
standard tableau. By identifying the words $w$ having the same tableau
$P^{A}(w)$, we obtain the plactic monoid $Pl(A_{n})$ whose defining relations
were determined by Knuth:
\begin{align*}
yzx  &  =yxz\text{ \ \ and \ \ }xzy=zxy\text{ if }x\prec y\prec z,\\
xyx  &  =xxy\text{ \ \ and \ \ }xyy=yxy\text{ if }x\prec y.
\end{align*}
The Robinson-Schensted correspondence has a natural interpretation in terms of
Kashiwara's theory of crystal bases \cite{DJM}, \cite{Ka2}, \cite{LLT}. Let
$V_{n}^{A}$ denote the vector representation of $U_{q}(sl_{n})$. By
considering each vertex of the crystal graph of $\underset{l\geq0}{\bigoplus
(}V_{n}^{A})^{\otimes l}$ as a word on $\mathcal{A}_{n}$, we have for any
words $w_{1}$ and $w_{2}$:

\begin{itemize}
\item $P^{A}(w_{1})=P^{A}(w_{2})$ if and only if $w_{1}$ and $w_{2}$ occur at
the same place in two isomorphic connected components of this graph.

\item $Q^{A}(w_{1})=Q^{A}(w_{2})$ if and only if $w_{1}$ and $w_{2}$ occur in
the same connected component of this graph.
\end{itemize}

Replacing $V_{n}^{A}$ by the vector representation $V_{n}^{C}$ of $sp_{2n}$
whose basis vectors are labelled by the letters of the totally ordered
alphabet
\[
\mathcal{C}_{n}=\{1\prec\cdot\cdot\cdot\prec n-1\prec n\prec\overline{n}%
\prec\overline{n-1}\prec\cdot\cdot\cdot\prec\overline{1}\},
\]
we have obtained in \cite{Lec} a Schensted type correspondence for type
$C_{n}$. This correspondence is based on an insertion algorithm for the
Kashiwara-Nakashima's symplectic tableaux \cite{KN} analogous to the bumping
algorithm. It may be regarded as a bijection between words $w$ of length $l$
on $\mathcal{C}_{n}$ and pairs $(P^{C}(w),Q^{C}(w))$ where $P^{C}(w)$ is a
symplectic tableau and $Q^{C}(w)$ an oscillating tableau of type $C$ and
length $l,$ that is, a sequence $(Q_{1},...,Q_{l})$ of Young diagrams such
that two consecutive diagrams differ by exactly one box. Moreover by
identifying the words of the free monoid $\mathcal{C}_{n}^{\ast}$ having the
same symplectic tableau we also obtain a monoid $Pl(C_{n})$. This is the
plactic monoid of type $C_{n}$ in the sense of \cite{Lit} and \cite{LLT}.

The vector representations $V_{n}^{B}$ and $V_{n}^{D}$ of $U_{q}(so_{2n+1})$
and $U_{q}(so_{2n})$ have crystal graphs whose vertices may be respectively
labelled by the letters of
\[
\mathcal{B}_{n}=\{1\prec\cdot\cdot\cdot\prec n-1\prec n\prec0\prec\overline
{n}\prec\overline{n-1}\prec\cdot\cdot\cdot\prec\overline{1}\}
\]
and
\[
\mathcal{D}_{n}=\{1\prec\cdot\cdot\cdot\prec n-1\prec%
\begin{array}
[c]{l}%
n\\
\overline{n}%
\end{array}
\prec\overline{n-1}\prec\cdot\cdot\cdot\prec\overline{1}\}.
\]
Let $G_{n}^{B}$ and $G_{n}^{D}$ be the crystal graphs of $\underset{l\geq
0}{\text{ }\bigoplus}(V_{n}^{B})^{\otimes l}$ and $\underset{l\geq0}%
{\bigoplus}(V_{n}^{D})^{\otimes l}.$ Then it is possible to label the vertices
of $G_{n}^{B}$ and $G_{n}^{D}$ by the words of the free monoids $\mathcal{B}%
_{n}^{\ast}$ and $\mathcal{D}_{n}^{\ast}$. However the situation is more
complicated than in the case of types $A$ and $C$. Indeed there exist a
fundamental representation of $U_{q}(so_{2n+1})$ and two fundamental
representations of $U_{q}(so_{2n})$ that do not appear in the decompositions
of $\bigoplus(V_{n}^{B})^{\otimes l}$ and $\underset{l\geq0}{\bigoplus}%
(V_{n}^{D})^{\otimes l}$ into their irreducible components. They are called
the spin representations and denoted respectively by $V(\Lambda_{n}^{B}),$
$V(\Lambda_{n}^{D})$ and $V(\Lambda_{n-1}^{D})$. In \cite{KN}, Kashiwara and
Nakashima have described their crystal graphs by using a new combinatorical
object that we will call a spin column. Write $SP_{n}$ for the set of spin
columns of height $n$ and set $\frak{B}_{n}=\mathcal{B}_{n}\cup SP_{n},$
$\frak{D}_{n}=\mathcal{D}_{n}\cup SP_{n}$. Then each vertex of the crystal
graphs $\frak{G}_{n}^{B}$ and $\frak{G}_{n}^{D}$ of $\underset{l\geq
0}{\bigoplus}\left(  V_{n}^{B}\bigoplus V(\Lambda_{n}^{B})\right)  ^{\otimes
l}$ and$\underset{l\geq0}{\text{ }\bigoplus}\left(  V_{n}^{D}\bigoplus
V(\Lambda_{n}^{D})\bigoplus V(\Lambda_{n-1}^{D})\right)  ^{\otimes l}$ may be
respectively identified with a word on $\frak{B}_{n}$ or $\frak{D}_{n}.$ We
can define two relations $\overset{B}{\sim}$ and $\overset{D}{\sim}$ by:

\begin{description}
\item $w_{1}\overset{B}{\sim}w_{2}$ if and only if $w_{1}$ and $w_{2}$ occur
at the same place in two isomorphic connected components of $\frak{G}_{n}^{B},$

\item $w_{1}\overset{D}{\sim}w_{2}$ if and only if $w_{1}$ and $w_{2}$ occur
at the same place in two isomorphic connected components of $\frak{G}_{n}^{D}.$
\end{description}

In this article, we prove that $Pl(B_{n})=\mathcal{B}_{n}^{\ast}/\overset
{B}{\sim},$ $Pl(D_{n})=\mathcal{D}_{n}^{\ast}/\overset{D}{\sim}$,
$\frak{Pl(}B_{n}\frak{)=B}_{n}^{\ast}\frak{/}\overset{B}{\sim}$ and
$\frak{Pl}(D_{n})=\frak{D}_{n}^{\ast}/\overset{D}{\sim}$ are monoids\ and we
undertake a detailed investigation of the corresponding insertion algorithms.
We summarize in part 2 the background on Kashiwara's theory of crystals used
in the sequel. In part 3, we first recall Kashiwara-Nakashima's notion of
orthogonal tableau (analogous to Young tableaux for types $B$ and $D$) and we
relate it to Littelmann's notion of Young tableau for classical types. Then we
derive a set of defining relations for $Pl(B_{n})$ and $Pl(D_{n})$ and we
describe the corresponding column insertion algorithms. Using the
combinatorial notion of oscillating tableaux (analogous to standard tableaux
for types $B$ and $D$), these algorithms yield the desired Schensted type
correspondences in $G_{n}^{B}$ and $G_{n}^{D}$. In part 4 we propose an
orthogonal Jeu de Taquin for type $B$ based on Sheats' sliding algorithm for
type $C$ \cite{SH}. Finally in part 5, we bring into the picture the spin
representations and extend the results of part 3 to the graphs $\frak{G}%
_{n}^{B}$, $\frak{G}_{n}^{D}$ and the monoids $\frak{Pl(}B_{n}\frak{)}$,
$\frak{Pl(}D_{n}\frak{).}$

\begin{notation}
In the sequel, we often write $B$ and $D$ instead of $B_{n}$ and $D_{n}$ to
simplify the notation. Moreover, we frequently define similar objects for
types $B$ and $D$. When they are related to type $B$ (respectively $D$), we
attach to them the label $^{B}$ (respectively the label $^{D}$). To avoid
cumbersome repetitions, we sometimes omit the labels $^{B}$ and $^{D}$ when
our statements are true for the two types.
\end{notation}

\section{Conventions for crystal graphs}

\subsection{Kashiwara's operators}

Let $\frak{g}$ be simple Lie algebra and $\alpha_{i},$ $i\in I$ its simple
roots. Recall that the crystal graphs of the $U_{q}(\frak{g})$-modules are
oriented colored graphs with colors $i\in I$. An arrow $a\overset
{i}{\rightarrow}b$ means that $\widetilde{f}_{i}(a)=b$ and $\widetilde{e}%
_{i}(b)=a$ where $\widetilde{e}_{i}$ and $\widetilde{f}_{i}$ are the crystal
graph operators (for a review of crystal bases and crystal graphs see
\cite{Ka2}). Let $V,V^{\prime}\ $be two $U_{q}(\frak{g})$-modules and
$B,B^{\prime}$ their crystal graphs. A vertex $v^{0}\in B$ satisfying
$\widetilde{e}_{i}(v^{0})=0$ for any $i\in I$ is called a highest weight
vertex. The decomposition of $V$ into its irreducible components is reflected
into the decomposition of $B$ into its connected components. Each connected
component of $B$ contains a unique vertex of highest weight. We write
$B(v^{0})$ for the connected component containing the highest weight vertex
$v^{0}$.\ The crystals graphs of two isomorphic irreducible components are
isomorphic as oriented colored graphs. We will say that two vertices $b_{1}%
\ $and $b_{2}$ of $B$ occur at the same place in two isomorphic connected
components $\Gamma_{1}$ and $\Gamma_{2}$ of $B$ if there exist $i_{1}%
,...,i_{r}\in I$ such that $w_{1}=\widetilde{f}_{i_{i}}\cdot\cdot
\cdot\widetilde{f}_{i_{r}}(w_{1}^{0})$ and $w_{2}=\widetilde{f}_{i_{i}}%
\cdot\cdot\cdot\widetilde{f}_{i_{r}}(w_{2}^{0})$, where $w_{1}^{0}$ and
$w_{2}^{0}$ are respectively the highest weight vertices of $\Gamma_{1}$ and
$\Gamma_{2}$.

The action of $\widetilde{e}_{i}$ and $\widetilde{f}_{i}$ on $B\otimes
B^{\prime}=\{b\otimes b^{\prime};$ $b\in B,b^{\prime}\in B^{\prime}\}$ is
given by:%

\begin{align}
\widetilde{f_{i}}(u\otimes v)  &  =\left\{
\begin{tabular}
[c]{c}%
$\widetilde{f}_{i}(u)\otimes v$ if $\varphi_{i}(u)>\varepsilon_{i}(v)$\\
$u\otimes\widetilde{f}_{i}(v)$ if $\varphi_{i}(u)\leq\varepsilon_{i}(v)$%
\end{tabular}
\right. \label{TENS1}\\
&  \text{and}\nonumber\\
\widetilde{e_{i}}(u\otimes v)  &  =\left\{
\begin{tabular}
[c]{c}%
$u\otimes\widetilde{e_{i}}(v)$ if $\varphi_{i}(u)<\varepsilon_{i}(v)$\\
$\widetilde{e_{i}}(u)\otimes v$ if $\varphi_{i}(u)\geq\varepsilon_{i}(v)$%
\end{tabular}
\right.  \label{TENS2}%
\end{align}
where $\varepsilon_{i}(u)=\max\{k;\widetilde{e}_{i}^{k}(u)\neq0\}$ and
$\varphi_{i}(u)=\max\{k;\widetilde{f}_{i}^{k}(u)\neq0\}$. Denote by
$\Lambda_{i},$ $i\in I$ the fundamental weights of $\frak{g}$. The weight of
the vertex $u$ is defined by $\mathrm{wt}(u)=\underset{I}{\sum}(\varphi
_{i}(u)-\varepsilon_{i}(u))\Lambda_{i}$. Write $s_{i}=s_{\alpha_{i}}$ for
$i\in I$.\ The Weyl group $W$ of $\frak{g}$ acts on $B$ by:
\begin{align}
s_{i}(u)  &  =(\widetilde{f_{i}})^{\varphi_{i}(u)-\varepsilon_{i}(u)}(u)\text{
if }\varphi_{i}(u)-\varepsilon_{i}(u)\geq0,\label{actionW}\\
s_{i}(u)  &  =(\widetilde{e_{i}})^{\varepsilon_{i}(u)-\varphi_{i}(u)}(u)\text{
if }\varphi_{i}(u)-\varepsilon_{i}(u)<0.\nonumber
\end{align}
We have the equality $\mathrm{wt}(\sigma(u))=\sigma(\mathrm{wt}(u)$ for any
$\sigma\in W$ and $u\in B.$ The following lemma is a straightforward
consequence of (\ref{TENS1}) and (\ref{TENS2}).

\begin{lemma}
\label{lem_phi_tens}Let $u\otimes v$ $\in$ $B\otimes B^{\prime}$. Then:

\begin{itemize}
\item $\mathrm{(i)}$ $\varphi_{i}(u\otimes v)=\left\{
\begin{tabular}
[c]{l}%
$\varphi_{i}(v)+\varphi_{i}(u)-\varepsilon_{i}(v)$ if $\varphi_{i}%
(u)>\varepsilon_{i}(v)$\\
$\varphi_{i}(v)$ otherwise.
\end{tabular}
\right.  .$

\item $\mathrm{(ii)}$ $\varepsilon_{i}(u\otimes v)=\left\{
\begin{tabular}
[c]{l}%
$\varepsilon_{i}(v)+\varepsilon_{i}(u)-\varphi_{i}(u)$ if $\varepsilon
_{i}(v)>\varphi_{i}(u)$\\
$\varepsilon_{i}(u)$ otherwise.
\end{tabular}
\right.  .$

\item $\mathrm{(iii)}$ $u\otimes v$ is a highest weight vertex of $B\otimes
B^{\prime}$ if and only if for any $i\in I$ $\widetilde{e}_{i}(u)=0$ (i.e. $u$
is of highest weight) and $\varepsilon_{i}(v)\leq\varphi_{i}(u).$
\end{itemize}
\end{lemma}

For any dominant weight $\lambda\in P_{+}$, write $B(\lambda)$\ for the
crystal graph of $V(\lambda),$ the irreducible module of highest weight
$\lambda$ and denote by $u_{\lambda}$ its highest weight vertex. Kashiwara has
introduced in \cite{Ka3} an embedding of $B(\lambda)$ into $B(m\lambda)$ for
any positive integer $m$. He uses this embedding to obtain a simple bijection
between Littelmann's path crystal associated to $\lambda$ and $B(\lambda)$
\cite{Lit3}.

\begin{theorem}
\label{th_strech}(Kashiwara) There exists a unique injective map
\begin{gather*}
S_{m}:B(\lambda)\rightarrow B(m\lambda)\subset B(\lambda)^{\otimes m}\\
u_{\lambda}\mapsto u_{\lambda}^{\otimes m}%
\end{gather*}
such that for any $b\in B(\lambda)$:
\begin{align}
\text{$\mathrm{(i)}$\ \ }S_{m}(\widetilde{e}_{i}(b))  &  =\widetilde{e}%
_{i}^{m}(S_{m}(b)),\nonumber\\
\text{$\mathrm{(ii)}$ \ }S_{m}(\widetilde{f}_{i}(b))  &  =\widetilde{f}%
_{i}^{m}(S_{m}(b)),\nonumber\\
\text{$\mathrm{(iii)}$ \ }\varphi_{i}(S_{m}(b))  &  =m\varphi_{i}%
(b),\label{stretch}\\
\text{$\mathrm{(iv)}$ \ }\varepsilon_{i}(S_{m}(b))  &  =m\varepsilon
_{i}(b),\nonumber\\
\text{$\mathrm{(v)}$ \ }\mathrm{wt}(S_{m}(b))  &  =m\mathrm{wt}(b).\nonumber
\end{align}
\end{theorem}

\begin{corollary}
\label{cor_strech}Let $\lambda_{1},...,\lambda_{k}\in P_{+}.$ Then, the map:
\begin{gather*}
\text{{\Large S}}_{m}:B(\lambda_{1})\otimes\cdot\cdot\cdot\otimes
B(\lambda_{k})\rightarrow B(m\lambda_{1})\otimes\cdot\cdot\cdot\otimes
B(m\lambda_{k})\\
\text{\ \ \ \ \ \ \ \ \ \ \ }b_{1}\otimes\cdot\cdot\cdot\otimes b_{k}\mapsto
S_{m}(b_{1})\otimes\cdot\cdot\cdot\otimes S_{m}(b_{k})
\end{gather*}
is injective and satisfies the relations (\ref{stretch}) with $b=b_{1}%
\otimes\cdot\cdot\cdot\otimes b_{k}.$ Moreover the image by {\Large S}$_{m}$
of a highest weight vertex of $B(\lambda_{1})\otimes\cdot\cdot\cdot\otimes
B(\lambda_{k})$ is a highest weight vertex of $B(m\lambda_{1})\otimes
\cdot\cdot\cdot\otimes B(m\lambda_{k})$.
\end{corollary}

\begin{proof}
By induction, we can suppose $k=2$. {\Large S}$_{m}$ is injective because
$S_{m}$ is injective. Let $u\otimes v\in B(\lambda_{1})\otimes B(\lambda
_{2}).$ Suppose that $\varphi_{i}(u)\leq\varepsilon_{i}(v)$. We derive the
following equalities from Formulas (\ref{TENS1}) and (\ref{TENS2}):
\begin{gather*}
\text{{\Large S}}_{m}\widetilde{f}_{i}(u\otimes v)=\text{{\Large S}}%
_{m}(u\otimes\widetilde{f}_{i}v)=S_{m}(u)\otimes S_{m}(\widetilde{f}%
_{i}v)=S_{m}(u)\otimes\widetilde{f}_{i}^{m}S_{m}(v)\\
\text{and }\widetilde{f}_{i}^{m}(\text{{\Large S}}_{m}(u\otimes v))=\widetilde
{f}_{i}^{m}(S_{m}(u)\otimes S_{m}(v))=S_{m}(u)\otimes\widetilde{f}_{i}%
^{m}S_{m}(v).
\end{gather*}
Indeed, $\varepsilon_{i}(S_{m}(v))=m\varepsilon_{i}(v)\geq m\varphi
_{i}(u)=\varphi_{i}(S_{m}(u))$ and for $p=1,...,m$ $\varepsilon_{i}%
(\widetilde{f}_{i}^{p}S_{m}(v))>\varepsilon_{i}(S_{m}(v)).$ Hence
{\Large S}$_{m}\widetilde{f}_{i}(u\otimes v)=\widetilde{f}_{i}^{m}($%
{\Large S}$_{m}(u\otimes v)).$ Now suppose $\varepsilon_{i}(v)<\varphi_{i}(u)$
i.e. $\varepsilon_{i}(u)\leq\varphi_{i}(v)+1$. We obtain:
\begin{gather*}
\text{{\Large S}}_{m}\widetilde{f}_{i}(u\otimes v)=\text{{\Large S}}%
_{m}(\widetilde{f}_{i}u\otimes v)=S_{m}(\widetilde{f}_{i}u)\otimes
S_{m}(v)=\widetilde{f}_{i}^{m}S_{m}(u)\otimes S_{m}(v)\\
\text{and }\widetilde{f}_{i}^{m}(\text{{\Large S}}_{m}(u\otimes v))=\widetilde
{f}_{i}^{m}(S_{m}(u)\otimes S_{m}(v))=\widetilde{f}_{i}^{m}S_{m}(u)\otimes
S_{m}(v)
\end{gather*}
because $\varepsilon_{i}(S_{m}(v))=m\varepsilon_{i}(v)\leq m\varphi
_{i}(u)+m=\varphi_{i}(S_{m}u)+m$. Hence we have {\Large S}$_{m}\widetilde
{f}_{i}(u\otimes v)=\widetilde{f}_{i}^{m}(${\Large S}$_{m}(u\otimes v)).$

Similarly we prove that {\Large S}$_{m}\widetilde{e}_{i}(u\otimes
v)=\widetilde{e}_{i}^{m}(${\Large S}$_{m}(u\otimes v)).$ So {\Large S}$_{m}$
satisfies the formulas $\mathrm{(i)}$\textrm{\ }and $\mathrm{(ii)}$. By Lemma
\ref{lem_phi_tens} $\mathrm{(i)}$\textrm{\ }and $\mathrm{(ii)}$ we obtain then
that {\Large S}$_{m}$ satisfies $\mathrm{(iii)}$, $\mathrm{(iv)}$ and
$\mathrm{(v)}$.

Suppose that $u\otimes v$ is a highest weight vertex of $B(\lambda_{1})\otimes
B(\lambda_{2})$.\ By Lemma \ref{lem_phi_tens} $\mathrm{(iii)}$, $u$ is the
highest weight vertex of $B(\lambda_{1})$ and $\varepsilon_{i}(v)\leq
\varphi_{i}(u)$ for $i\in I$.\ Then by definition of $S_{m},$ $S_{m}(u)$ is
the highest weight vertex of $B(m\lambda_{1})$.\ Moreover for any $i\in I,$
$\varepsilon_{i}(S_{m}(v))=m\varepsilon_{i}(v)\leq m\varphi_{i}(u)=\varphi
_{i}(S_{m}(u))$. So $S_{m}(u)\otimes S_{m}(v)=${\Large S}$_{m}(u\otimes v)$ is
of highest weight in $B(m\lambda_{1})\otimes B(m\lambda_{2})$.
\end{proof}

\bigskip

\noindent By this corollary, the connected component of $B(\lambda_{1}%
)\otimes\cdot\cdot\cdot\otimes B(\lambda_{k})$ of highest weight vertex
$u^{0}=u_{1}\otimes\cdot\cdot\cdot\otimes u_{k}$, may be identified with the
sub-graph of $B(m\lambda_{1})\otimes\cdot\cdot\cdot\otimes B(m\lambda_{k})$
generated by the vertex $S_{m}(u_{1})\otimes\cdot\cdot\cdot\otimes S_{m}%
(u_{k})$ and the operators $\widetilde{f}_{i}^{m}$ for $i\in I$.

\subsection{Tensor powers of the vector representations}

We choose to label the Dynkin diagram of $so_{2n+1}$ by:
\[
\overset{1}{\circ}-\overset{2}{\circ}-\overset{3}{\circ}\cdot\cdot
\cdot\overset{n-2}{\circ}-\overset{n-1}{\circ}\Longrightarrow\overset{n}%
{\circ}%
\]
and the Dynkin diagram of $so_{2n}$ by:%

\[
\overset{1}{\circ}-\overset{2}{\circ}-\overset{3}{\circ}\cdot\cdot
\cdot\overset{n-3}{\circ}-
\begin{tabular}
[c]{l}%
$\ \ \ \ \overset{n}{\circ}$\\
$\ \ \ /$\\
$\overset{n-2}{\circ}$\\
$\ \ \ \backslash$\\
\ \ \ $\underset{n-1}{\circ}$%
\end{tabular}
.
\]
Write $W_{n}^{B}$ and $W_{n}^{D}$ for the Weyl groups of $so_{2n+1}$ and
$so_{2n}$. Denote by $V_{n}^{B}$ and $V_{n}^{D}$ the vector representations of
$U_{q}(so_{2n+1})$ and $U_{q}(so_{2n}).$ Their crystal graphs are
respectively:
\begin{equation}
1\overset{1}{\rightarrow}2\cdot\cdot\cdot\rightarrow n-1\overset
{n-1}{\rightarrow}n\overset{n}{\rightarrow}0\overset{n}{\rightarrow}%
\overline{n}\overset{n-1}{\rightarrow}\overline{n-1}\overset{n-2}{\rightarrow
}\cdot\cdot\cdot\rightarrow\overline{2}\overset{1}{\rightarrow}\overline{1}
\label{vect_B}%
\end{equation}
and
\begin{equation}
1\overset{1}{\rightarrow}2\overset{2}{\rightarrow}\cdot\cdot\cdot\overset
{n-3}{\rightarrow}n-2\overset{n-2}{\rightarrow}
\begin{tabular}
[c]{c}%
$\overline{n}$ \ \ \\
\ \ $\overset{n}{\nearrow}$ $\ \ \ \overset{n-1}{\text{ \ }\searrow}$ \ \ \ \\
$n-1\ \ \ \ \ \ \ \ \ \ \overline{n-1}$\\
\ $\underset{n-1}{\searrow}$ \ \ \ $\underset{n}{\nearrow}$ \ \ \ \\
$n$ \
\end{tabular}
\overset{n-2}{\rightarrow}\overline{n-2}\overset{n-3}{\rightarrow}\cdot
\cdot\cdot\overset{2}{\rightarrow}\overline{2}\overset{1}{\rightarrow
}\overline{1}. \label{vect_D}%
\end{equation}
By induction, formulas (\ref{TENS1}), (\ref{TENS2}) allow to define a crystal
graph for the representations $(V_{n}^{B})^{\otimes l}$ and $(V_{n}%
^{D})^{\otimes l}$ for any $l$. Each vertex $u_{1}\otimes u_{2}\otimes
\cdot\cdot\cdot\otimes u_{l}$ of the crystal graph of $(V_{n}^{B})^{\otimes
l}$ will be identified with the word $u_{1}u_{2}\cdot\cdot\cdot u_{l}$ on the
totally ordered alphabet
\[
\mathcal{B}_{n}=\{1\prec\cdot\cdot\cdot\prec n-1\prec n\prec0\prec\overline
{n}\prec\overline{n-1}\prec\cdot\cdot\cdot\prec\overline{1}\}.
\]
Similarly each vertex $v_{1}\otimes v_{2}\otimes\cdot\cdot\cdot\otimes v_{l}$
of the crystal graph of $(V_{n}^{D})^{\otimes l}$ will be identified with the
word $v_{1}v_{2}\cdot\cdot\cdot v_{l}$ on the partially ordered alphabet
\[
\mathcal{D}_{n}=\{1\prec\cdot\cdot\cdot\prec n-1\prec%
\begin{array}
[c]{l}%
n\\
\overline{n}%
\end{array}
\prec\overline{n-1}\prec\cdot\cdot\cdot\prec\overline{1}\}.
\]
By convention we set $\overline{0}=0$ and for $k=1,\cdot\cdot\cdot,n,$
$\overline{\overline{k}}=k$. The letter $x$ is barred if $x\succeq\overline
{n}$ unbarred if $x\preceq n$ and we set:
\[
\left|  x\right|  =\left\{
\begin{tabular}
[c]{l}%
$x$ if $x$ is unbarred\\
$\overline{x}$ otherwise.
\end{tabular}
\right.
\]
Write $\mathcal{B}_{n}^{\ast}$ and $\mathcal{D}_{n}^{\ast}$ for the free
monoids on $\mathcal{B}_{n}$ and $\mathcal{D}_{n}$. If $w$ is a word of
$\mathcal{B}_{n}^{\ast}$ or $\mathcal{D}_{n}^{\ast}$, we denote by
$\mathrm{l}(w)$ its length and by $d(w)=(d_{1},...,d_{n})$ the $n$-tuple where
$d_{i}$ is the number of letters $i$ in $w$ minus the number of letters
$\overline{i} $. Let $G_{n}^{B}$ and $G_{n,l}^{B}$ be respectively the crystal
graphs of $\underset{l}{\bigoplus}(V_{n}^{B})^{\otimes l}$ and $(V_{n}%
^{B})^{\otimes l}$. Then the vertices of $G_{n}^{B}$ are indexed by the words
of $\mathcal{B}_{n}^{\ast}$ and those of $G_{n,l}^{B}$ by the words of
$\mathcal{B}_{n}^{\ast}$ of length $l$. Similarly $G_{n}^{D}$ and $G_{n,l}%
^{D}$, the crystal graphs of $\underset{l}{\bigoplus}(V_{n}^{B})^{\otimes l}$
and $(V_{n}^{B})^{\otimes l}$ are indexed respectively by the words of
$\mathcal{D}_{n}^{\ast}$ and by the words of $\mathcal{D}_{n}^{\ast}$ of
length $l$. If $w$ is a vertex of $G_{n}$, write $B(w)$ for the connected
component of $G_{n}$ containing $w$.

Denote by $\Lambda_{1}^{B},...,\Lambda_{n}^{B}$ and $\Lambda_{1}%
^{D}...,\Lambda_{n}^{D}$ the fundamental weights of $U_{q}(so_{2n+1})$ and
$U_{q}(so_{2n}).$ Let $P_{+}^{B}$ and $P_{+}^{D}$ be the sets of dominant
weights of their weight lattices. We set
\begin{align*}
\omega_{n}^{B}  &  =2\Lambda_{n}^{B},\\
\omega_{i}^{B}  &  =\Lambda_{i}^{B}\text{ for }i=1,...,n-1
\end{align*}
and
\begin{align*}
\omega_{n}^{D}  &  =2\Lambda_{n}^{D},\\
\overline{\omega}_{n}^{D}  &  =2\Lambda_{n-1}^{D},\\
\omega_{n-1}^{D}  &  =\Lambda_{n}^{D}+\Lambda_{n-1}^{D},\\
\omega_{i}^{D}  &  =\Lambda_{i}^{D}\text{ for }i=1,...,n-2.
\end{align*}
Then the weight of a vertex $w$ of $G_{n}$ is given by:
\[
\mathrm{wt}(w)=d_{n}\omega_{n}+\overset{n-1}{\underset{i=1}{\sum}}%
(d_{i}-d_{i+1})\omega_{i}.
\]
Thus we recover the well-known fact that there is no connected component of
$G_{n}^{B}$ isomorphic to $B(\Lambda_{n}^{B})$ and no connected component of
$G_{n}^{D}$ isomorphic to $B(\Lambda_{n}^{D})$ or $B(\Lambda_{n-1}^{D})$.
Recall that in the cases of the types $A$ and $C,$ every crystal graph of an
irreducible module may be embedded in the crystal graph of a tensor power of
the vector representation. For $\lambda\in P_{+}^{B}$,$\ B^{B}(\lambda)$ may
be embedded in a tensor power of the vector representation $V_{n}^{B}$ if and
only if $\lambda$ lies in the weight sub-lattice $\Omega^{B}$ generated by the
$\omega_{i}^{B}$'s. Similarly, for $\lambda\in P_{+}^{D},$ $B^{D}(\lambda)$
may be embedded in a tensor power of the vector representation $V_{n}^{D}$ if
and only if $\lambda$ lies in the weight sub-lattice $\Omega^{D}$ generated by
the $\omega_{i}^{D}$'s.\ Set $\Omega_{+}^{B}=P_{+}^{B}\cap\Omega^{B}$ and
$\Omega_{+}^{D}=P_{+}^{D}\cap\Omega^{D}$.

Now we introduce the coplactic relation. For $w_{1}$ and $w_{2}\in$
$\mathcal{B}_{n}^{\ast}$ (resp. $\mathcal{D}_{n}^{\ast}$), write
$w_{1}\overset{B}{\longleftrightarrow}w_{2}\ $(resp. $w_{1}\overset
{D}{\longleftrightarrow}w_{2}$) if and only if $w_{1}$ and $w_{2}$ belong to
the same connected component of $G_{n}^{B}$ (resp. $G_{n}^{D}$).\ The proof of
the following lemma is the same as in the symplectic case \cite{Lec}.

\begin{lemma}
\label{lem_coplactic}If $w_{1}=u_{1}v_{1}$ and $w_{2}=u_{2}v_{2}$ with
$l(u_{1})=l(u_{2})$ and $l(v_{1})=l(v_{2})$%
\[
w_{1}\longleftrightarrow w_{2}\Longrightarrow\left\{
\begin{tabular}
[c]{l}%
$u_{1}\longleftrightarrow u_{2}$\\
$v_{1}\longleftrightarrow v_{2}$%
\end{tabular}
\right.  .
\]
\end{lemma}

\subsection{Crystal graphs of the spin representations}

The spin representations of $U_{q}(so_{2n+1})$ and $U_{q}(so_{2n})$ are
$V(\Lambda_{n}^{B}),$ $V(\Lambda_{n}^{D})$ and $V(\Lambda_{n-1}^{D})$. Recall
that $\dim V(\Lambda_{n}^{B})=2^{n}$ and $\dim V(\Lambda_{n}^{D})=\dim
V(\Lambda_{n-1}^{D})=2^{n-1}$. Now we review the description of $B(\Lambda
_{n}^{B}),$ $B(\Lambda_{n}^{D})$ and $B(\Lambda_{n-1}^{D})$ given by Kashiwara
and Nakashima in \cite{KN}. It is based on the notion of spin column. To avoid
confusion between these new columns and the classical columns of a tableau
that we introduce in the next section, we follow Kashiwara-Nakashima's
convention and represent spin columns by column shape diagrams of width $1/2$.
Such diagrams will be called K-N diagrams.

\begin{definition}
A spin column $\frak{C}$ of height $n$ is a K-N diagram containing $n$ letters
of $\mathcal{D}_{n}$ such that the word $x_{1}\cdot\cdot\cdot x_{n}$ obtained
by reading $\frak{C}$ from top to bottom does not contain a pair
$(z,\overline{z})$ and verifies $x_{1}\prec\cdot\cdot\cdot\prec x_{n}$. The
set of spin columns of length $n$ will be denoted $SP_{n}$.
\end{definition}

\begin{itemize}
\item $B(\Lambda_{n}^{B})=\{\frak{C;}$ $\frak{C}\in SP_{n}\}$ where
Kashiwara's operators act as follows:

\begin{description}
\item  if $n\in\frak{C}$ then $\widetilde{f}_{n}\frak{C}$ is obtained by
turning $n$ into $\overline{n}$, otherwise $\widetilde{f}_{n}\frak{C}=0$,

\item  if $\overline{n}\in\frak{C}$ then $\widetilde{e}_{n}\frak{C}$ is
obtained by turning $\ \overline{n}$ into $n,$ otherwise $\widetilde{e}%
_{n}\frak{C}=0$,

\item  if $(i,\overline{i+1})\in\frak{C}$ then $\widetilde{f}_{i}\frak{C}$ is
obtained by turning $(i,\overline{i+1})$ into $(i+1,\overline{i})$, otherwise
$\widetilde{f}_{i}\frak{C}=0$,

\item  if $(i+1,\overline{i})\in\frak{C}$ then $\widetilde{e}_{i}\frak{C}$ is
obtained by turning $(i+1,\overline{i})$ into $(i,\overline{i+1})$, otherwise
$\widetilde{e}_{i}\frak{C}=0$.
\end{description}

\item $B(\Lambda_{n}^{D})=\{\frak{C}\in SP_{n};$ the number of barred letters
in $\frak{C}$ is even$\}$ and $B(\Lambda_{n-1}^{D})=\{\frak{C}\in SP_{n};$ the
number of barred letters in $\frak{C}$ is odd$\}$ where Kashiwara's operators
act as follows:

\begin{description}
\item  if $(n-1,n)\in\frak{C}$ then $\widetilde{f}_{n}\frak{C}$ is obtained by
turning $(n-1,n)\ $into $(\overline{n},\overline{n-1})$, otherwise
$\widetilde{f}_{n}\frak{C}=0$,

\item  if $(\overline{n},\overline{n-1})\in\frak{C}$ then $\widetilde{e}%
_{n}\frak{C}$ is obtained by turning $(\overline{n},\overline{n-1})$ into
$(n-1,n)$, otherwise $\widetilde{e}_{n}\frak{C}=0$,

\item  for $i\neq n,$ $\widetilde{f}_{i}$ and $\widetilde{e}_{i}$ act like in
$B(\Lambda_{n}^{B}).$
\end{description}
\end{itemize}%

\begin{figure}
[ptb]
\begin{center}
\includegraphics[
height=7.3301cm,
width=10.1067cm
]%
{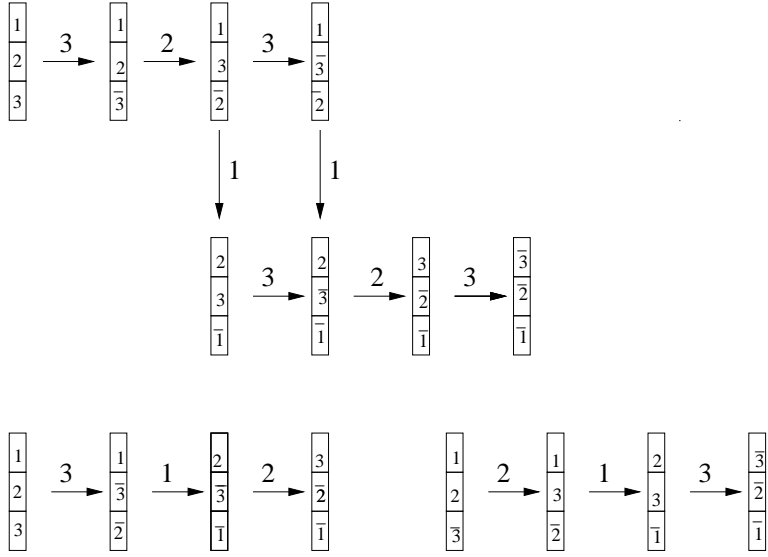}%
\caption{The crystal graphs $B(\Lambda_{n}^{B}),$ $B(\Lambda_{n}^{D})$ and
$B(\Lambda_{n-1}^{D})$ for $U_{q}(so_{7})$ and $U_{q}(so_{6})$}%
\end{center}
\end{figure}

In the sequel we denote by $v_{\Lambda_{n}}^{B}$ the highest weight vertex of
$B(\Lambda_{n}^{B})$, by $v_{\Lambda_{n}}^{D}$ and $v_{\Lambda_{n-1}}^{D}$ the
highest weight vertices of $B(\Lambda_{n}^{D})$ and $B(\Lambda_{n-1}^{D}) $.
Note that $v_{\Lambda_{n}}^{B}$ and $v_{\Lambda_{n}}^{D}$ correspond to the
spin column containing the letters of $\{1,...,n\}$ and $v_{\Lambda_{n-1}}%
^{D}$ corresponds to the spin column containing the letters of
$\{1,...,n-1,\overline{n}\}$.

\section{\label{sec_in_vect}Schensted correspondences in $G_{n}^{B}$ and
$G_{n}^{D}$}

\subsection{Orthogonal tableaux}

Let $\lambda\in\Omega_{+}$. We are going to review the notion of standard
orthogonal tableaux introduced by Kashiwara and Nakashima \cite{KN} to label
the vertices of $B(\lambda)$.

\subsubsection{\label{subsubsec_colo}Columns and admissible columns}

A column of type $B$ is a Young diagram
\[
C=
\begin{tabular}
[c]{|l|}\hline
$x_{1}$\\\hline
$\cdot$\\\hline
$\cdot$\\\hline
$x_{l}$\\\hline
\end{tabular}
\]
of column shape filled by letters of $\mathcal{B}_{n}$ such that $C$ increases
from top to bottom and $0$ is the unique letter of $\mathcal{B}_{n} $ that may
appear more than once.

\noindent A column of type $D$ is a Young diagram $C$ of column shape filled
by letters of $\mathcal{D}_{n}$ such that $x_{i+1}\nleqslant x_{i}$ for
$i=1,...,l-1$. Note that the letters $n$ and $\overline{n}$ are the unique
letters that may appear more than once in $C$ and if they do, these letters
are different in two adjacent boxes.

\noindent The height $h(C)$ of the column $C$ is the number of its letters.
The word obtained by reading the letters of $C$ from top to bottom is called
the reading of $C$ and denoted by \textrm{w}$(C)$. We will say that the column
$C$ contains a pair $(z,\overline{z})$ when a letter $0$ or the two letters
$z\preceq n$ and $\overline{z}$ appear in $C$.

\begin{definition}
(Kashiwara-Nakashima) Let $C$ be a column such that \textrm{w}$(C)=x_{1}%
\cdot\cdot\cdot x_{h(C)}$.\ Then $C$ is admissible if $h(C)\leq n$ and for any
pair $(z,\overline{z})$ of letters in $C$ satisfying $z=x_{p}$ and
$\overline{z}=x_{q}$ with $z\preceq n$ we have
\begin{equation}
\left|  q-p\right|  \geq h(C)-z+1. \label{cond_admissi}%
\end{equation}
(Note that $0\succ n$ on $\mathcal{B}_{n}$ and we may have $q-p<0$ for type
$D$ and $z=n).$
\end{definition}

\begin{example}
For $n=4,$ $40\bar{4}\bar{2}$ and $3\overline{4}4\overline{3}$ are readings of
admissible columns respectively of type $B$ and $D$.
\end{example}

\noindent Let $C$ be a column of type $B$ or $D$ and $z\preceq n$ a letter of
$C$. We denote by $N(z)$ the number of letters $x$ in $C$ such that $x\preceq
z$ or $x\succeq\overline{z}$. Then Condition (\ref{cond_admissi}) is
equivalent to $N(z)\leq z.$

\noindent Suppose that $C$ is non admissible and does not contain a pair
$(z,\overline{z})$ with $z\preceq n$ and $N(z)>z.$ Then $h(C)>n$. Hence $C$ is
of type $B$ and $0\in C.$ Indeed, if $0\notin C,$ $C$ contains a letter $z$
maximal such that $z\preceq n$ and $\overline{z}\in C$. It means that for any
$x\in\{z+1,...,n\},$ there is at most one letter $y\in C$ with $\left|
y\right|  =x.$ We have a contradiction because in this case $N(z)>n-(n-z)$. We
obtain the

\begin{remark}
\label{not_N(z)}A column $C$ is non admissible if and only if at least one of
the following assertions is satisfied:

$\mathrm{(i)}:C$ contains a letter $z\preceq n$ and $N(z)>z$

$\mathrm{(ii)}:C$ is of type $B,$ $0\in C$ and $h(C)>n.$
\end{remark}

\noindent If we set $v_{\omega_{k}}^{B}=1\cdot\cdot\cdot k$ for $k=1,...,n,$
then $B(v_{\omega_{k}}^{B})$ is isomorphic to $B(\omega_{k}^{B})$. Similarly,
if we set $v_{\omega_{k}}^{D}=1\cdot\cdot\cdot k$ for $k=1,...,n$ and
$v_{\overline{\omega}_{n}}^{D}=1\cdot\cdot\cdot(n-1)\overline{n},$ then
$B(v_{\omega_{k}}^{D})$ and $B(v_{\overline{\omega}_{n}}^{D})$ are
respectively isomorphic to $B(\omega_{k}^{D})$ and $B(\overline{\omega}%
_{n}^{D})$.

\begin{proposition}
\label{prop_KNf}(Kashiwara-Nakashima)

\begin{itemize}
\item  The vertices of $B(v_{\omega_{k}}^{B})$ are the readings of the
admissible columns of type $B$ and length $k$.

\item  The vertices of $B(v_{\omega_{k}}^{D})$ with $k<n$ are the readings of
the admissible columns of type $D$ and length $k$.

\item  The vertices of $B(v_{\omega_{n}}^{D})$ are the readings of the
admissible columns$\ C$ of type $D$ such that $\mathrm{w}(C)=x_{1}\cdot
\cdot\cdot x_{n}$ and $x_{k}=n$ (resp.\ $x_{k}=\overline{n}$) implies $n-k$ is
even (resp. odd).

\item  The vertices of $B(v_{\overline{\omega}_{n}}^{D})$ are the readings of
the admissible columns $C$ of type $D$ such that $\mathrm{w}(C)=x_{1}%
\cdot\cdot\cdot x_{n}$ and $x_{k}=\overline{n}$ (resp. $x_{k}=n$) implies
$n-k$ is odd (resp. even).
\end{itemize}
\end{proposition}

We can obtain another description of the admissible columns by computing, for
each admissible column $C$, a pair of columns $(lC,rC)$ without pair
$(z,\overline{z})$. This duplication was inspired by the description of the
admissible columns of type $C$ in terms of De Concini columns used by Sheats
in \cite{SH}.

\begin{definition}
Let $C$ be a column of type $B$ and denote by $I_{C}=\{z_{1}=0,...,z_{r}%
=0\succ z_{r+1}\succ\cdot\cdot\cdot\succ z_{s}\}$ the set of letters
$z\preceq0$ such that the pair $(z,\overline{z})$ occurs in $C$. We will say
that $C$ can be split when there exists (see the example below) a set of $s$
unbarred letters $J_{C}=\{t_{1}\succ\cdot\cdot\cdot\succ t_{s}\}\subset
\mathcal{B}_{n}$ such that:

$\ \ \ t_{1}$ is the greatest letter of $\mathcal{B}_{n}$ satisfying:
$t_{1}\prec z_{1},t_{1}\notin C$ and $\overline{t_{1}}\notin C,$

\ \ \ for $i=2,...,s$, $t_{i}$ is the greatest letter of $\mathcal{B}_{n}$
satisfying: $t_{i}\prec\min(t_{i-1,}z_{i}),$ $t_{i}\notin C$ and
$\overline{t_{i}}\notin C.$

\noindent In this case we write:

\begin{itemize}
\item $rC$ for the column obtained first by changing in $C$ $\overline{z}_{i}$
into $\overline{t}_{i}$ for each letter $z_{i}\in I,$ next by reordering if necessary.

\item $lC$ for the column obtained first by changing in $C$ $z_{i}$ into
$t_{i}$ for each letter $z_{i}\in I,$ next by reordering if necessary.
\end{itemize}
\end{definition}

\begin{definition}
\label{def_C_hat}Let $C$ be a column of type $D$. Denote by $\widehat{C}$ the
column of type $B$ obtained by turning in $C$ each factor $\overline{n}n$ into
$00$.\ We will say that $C$ can be split when $\widehat{C}$ can be split. In
this case we write $lC=l\widehat{C}$ and $rC=l\widehat{C}.$
\end{definition}

\begin{example}
\label{exam_splitting}Suppose $n=9$ and consider the column $C$ of type $B$
such that $\mathrm{w}(C)=458900\bar{8}\bar{5}\bar{4}$. We have $I_{C}%
=\{0,0,8,5,4\}$ and $J_{C}=\{7,6,3,2,1\}$. Hence
\[
\mathrm{w}(lC)=123679\bar{8}\bar{5}\bar{4}\text{ and }\mathrm{w}%
(rC)=4589\bar{7}\bar{6}\bar{3}\bar{2}\bar{1}.
\]
Suppose $n=8$ and consider the column $C^{\prime}$ of type $D$ such that
$\mathrm{w}(C^{\prime})=56\bar{8}8\bar{8}\bar{6}\bar{5}\bar{2}$.\ Then
$\mathrm{w}(\widehat{C^{\prime}})=5600\bar{8}\bar{6}\bar{5}\bar{2},$
$I_{\widehat{C^{\prime}}}=\{0,0,6,5\}$ and $J_{\widehat{C^{\prime}}%
}=\{7,4,3,1\}$. Hence
\[
\mathrm{w(}lC^{\prime})=1347\bar{8}\bar{6}\bar{5}\bar{2}\text{ and
}\mathrm{w(}rC^{\prime})=56\bar{8}\bar{7}\bar{4}\bar{3}\bar{2}\bar{1}.
\]
\end{example}

\begin{lemma}
\label{Lem_dC_impl_Cadm}Let $C$ be a column of type $B$ or $D$ which can be
split . Then $C$ is admissible.
\end{lemma}

\begin{proof}
Suppose $C$ of type $B$. We have $h(C)\leq n$ for $C$ can be split. If there
exists a letter$\ z\prec0$ in $C$ such that the pair $(z,\overline{z})$ occurs
in $C$ and $N(z)\geq z+1$, $C$ contains at least $z+1$ letters $x$ satisfying
$\left|  x\right|  \preceq z$. So $rC$ contains at least $z+1$ letters
$x^{\prime}$ satisfying $\left|  x^{\prime}\right|  \preceq z$. We obtain a
contradiction because $rC$ does not contain a pair $(t,\overline{t})$. When
$C$ is of type $D,$ by applying the lemma to $\widehat{C}$ we obtain that
$\widehat{C}$ is admissible. So $C$ is admissible.
\end{proof}

The meaning of $lC$ and $rC$ is explained in the following proposition.

\begin{proposition}
\label{prop_imag_S2}Let $\omega\in\{\omega_{1}^{B},...,\omega_{n}^{B})$ or
$\omega\in\{\omega_{1}^{D},...,\omega_{n-1}^{D},\omega_{n}^{D},\overline
{\omega}_{n}^{D}\}$. The map
\[
S_{2}:B(v_{\omega})\rightarrow B(v_{\omega})\otimes B(v_{\omega})
\]
defined in Theorem \ref{th_strech} satisfies for any admissible column $C\in
B(v_{\omega})$:
\[
S_{2}(\mathrm{w(}C))=\mathrm{w(}rC)\otimes\mathrm{w(}lC).
\]
\end{proposition}

\begin{example}
Consider $\omega=\omega_{2}^{B}$ for $U_{q}(so_{5})$. The following graphs are
respectively those of $B(\omega)$ and $S_{2}(B(\omega)).$%
\[%
\begin{tabular}
[c]{l}%
$12\overset{2}{\rightarrow}10\overset{2}{\rightarrow}1\bar{2}\overset
{1}{\rightarrow}2\bar{2}\overset{1}{\rightarrow}2\bar{1}$\\
\ \ \ \ \ \ \ \ $\downarrow1$ \ \ \ \ \ \ \ \ \ \ \ \ \ \ \ \ $\downarrow2$\\
$\text{\ \ \ \ \ \ \ }20\overset{2}{\rightarrow}00\overset{2}{\rightarrow
}0\bar{2}\overset{1}{\rightarrow}0\bar{1}\overset{2}{\rightarrow}\bar{2}%
\bar{1}$%
\end{tabular}
\]%
\[%
\begin{tabular}
[c]{l}%
$(12)\otimes(12)\overset{2^{2}}{\rightarrow}(1\bar{2})\otimes(12)\overset
{2^{2}}{\rightarrow}(1\bar{2})\otimes(1\bar{2})\overset{1^{2}}{\rightarrow
}(2\bar{1})\otimes(1\bar{2})\overset{1^{2}}{\rightarrow}(2\bar{1}%
)\otimes(2\bar{1})$\\
\ \ \ \ \ \ \ \ \ \ \ \ \ \ \ \ \ \ \ \ \ \ \ \ $\downarrow1^{2}$
\ \ \ \ \ \ \ \ \ \ \ \ \ \ \ \ \ \ \ \ \ \ \ \ \ \ \ \ \ \ \ \ \ \ \ \ \ \ \ \ \ \ \ \ \ $\downarrow
2^{2}$\\
$\ \ \ \ \ \ \ \ \ \ \ \ \ \ \ \ \ (2\bar{1})\otimes(12)\overset{2^{2}%
}{\rightarrow}(\bar{2}\bar{1})\otimes(12)\overset{2^{2}}{\rightarrow}(\bar
{2}\bar{1})\otimes(1\bar{2})\overset{1^{2}}{\rightarrow}(\bar{2}\bar
{1})\otimes(2\bar{1})\overset{2^{2}}{\rightarrow}(\bar{2}\bar{1})\otimes
(\bar{2}\bar{1})$%
\end{tabular}
\]
\end{example}

\begin{proof}
(of proposition \ref{prop_imag_S2}) In this proof we identify each column with
its reading to simplify the notations. When $C=v_{\omega}$ is the highest
weight vertex of $B(v_{\omega})$, $r(v_{\omega})=l(v_{\omega})=v_{\omega}$
because $v_{\omega}$ does not contain a pair $(z,\overline{z})$. So
$S_{2}(v_{\omega})=rC\otimes lC.$ Each vertex $C$ of $B(\omega)$ may be
written $C=\widetilde{f}_{i_{1}}\cdot\cdot\cdot\widetilde{f}_{i_{r}}%
(v_{\omega})$. By induction on $r$, it suffices to prove that for any
$\mathrm{w}(C)\in B(v_{\omega})$ such that $\widetilde{f}_{i}(C)\neq0$ we
have
\[
S_{2}(C)=rC\otimes lC\Longrightarrow S_{2}(\widetilde{f}_{i}C)=r(\widetilde
{f}_{i}C)\otimes l(\widetilde{f}_{i}C).
\]
For any column $D$ we denote by $[D]_{i}$ the word obtained by erasing all the
letters $x$ of $D$ such that $\widetilde{f}_{i}(x)=\widetilde{e}_{i}(x)=0$. It
is clear that only the letters of $[D]_{i}$ may be changed in $D$ when we
apply $\widetilde{f}_{i}$.

\noindent Suppose $\omega\in\{\omega_{1}^{B},...,\omega_{n}^{B})$. Consider
$C\in B(v_{\omega})$ such that $S_{2}(C)=rC\otimes lC$ and $\widetilde{f}%
_{i}(C)\neq0$.

\noindent When $i\neq n$, the letters $x\notin\{\overline{i+1},\overline
{i},i,i+1\}$ do not interfere in the computation of $\widetilde{f}_{i}$. It
follows from the condition $\widetilde{f}_{i}(C)\neq0$ and an easy computation
from (\ref{TENS1}) and (\ref{TENS2}) that we need only consider the following
cases: $\mathrm{(i)}$\textrm{\ }$[C]_{i}=i$, $\mathrm{(ii)}$ $[C]_{i}%
=\overline{i+1}$, $\mathrm{(iii)}$ $[C]_{i}=(i+1)(\overline{i+1}),$
$\mathrm{(iv)}$ $[C]_{i}=(i)(\overline{i+1}),$ $\mathrm{(v)}$ $[C]_{i}%
=i(i+1)(\overline{i+1})$ and $\mathrm{(vi)}$ $[C]_{i}=i(\overline
{i+1})\overline{i}$. In the case $\mathrm{(i)}$, if $i+1\notin J_{C},$ we have
$[lC]_{i}=i$ and $[rC]_{i}=i$.\ Then $[\widetilde{f}_{i}(C)]_{i}=i+1$ and
$J_{\widetilde{f}_{i}C}=J_{C}$ (hence $i\notin J_{\widetilde{f}_{i}C}$). So
$[l(\widetilde{f}_{i}C)]_{i}=i+1$ and $[r(\widetilde{f}_{i}C)]_{i}=i+1.$ That
means that $S_{2}(\widetilde{f}_{i}C)=\widetilde{f}_{i}^{2}(rC\otimes
lC)=\widetilde{f}_{i}(rC)\otimes\widetilde{f}_{i}(lC)=r(\widetilde{f}%
_{i}C)\otimes l(\widetilde{f}_{i}C)$ by definition of the map $S_{2}$. If
$i+1\in J_{C}$, we can write $[rC]_{i}=(i)(\overline{i+1})$ and $[lC]_{i}%
=(i)(i+1)$. Then $[\widetilde{f}_{i}C)]_{i}=i+1$ and $J_{\widetilde{f}_{i}%
C}=J_{C}-\{i+1\}+\{i\}$. So $[r(\widetilde{f}_{i}C)]=(i+1)(\overline{i})$ and
$[l(\widetilde{f}_{i}C)]=(i)(i+1).$ Hence $S_{2}(\widetilde{f}_{i}%
C)=\widetilde{f}_{i}^{2}(rC\otimes lC)=\widetilde{f}_{i}^{2}(rC)\otimes
lC=r(\widetilde{f}_{i}C)\otimes l(\widetilde{f}_{i}C)$. The proof is similar
in the cases $\mathrm{(ii)}$ to $\mathrm{(vi)}$. When $i=n$, only the letters
of $\{\overline{n},0,n\}$ interfere in the computation of $\widetilde{f}_{n}$.
We obtain the proposition by considering the cases: $[C]_{n}=\underset{0\text{
p times}}{\underbrace{0\cdot\cdot\cdot0}},$ $[C]_{n}=n\underset{0\text{ p
times}}{\underbrace{0\cdot\cdot\cdot0}}$ and $[C]_{n}=n$.

Suppose $\omega\in\{\omega_{1}^{D},...,\omega_{n-1}^{D},\overline{\omega}%
_{n}^{D},\omega_{n}^{D}\}$. When $i<n-1$ the proof is the same than above.
When $i\in\{n-1,n\},$ the proposition follows by considering successively the
cases:
\[
\left\{
\begin{tabular}
[c]{l}%
$\lbrack C]_{i}=n-1(\overline{n}n)^{r}$,\\
$\lbrack C]_{i}=n(\overline{n}n)^{r}\overline{n}$,\\
$\lbrack C]_{i}=(n-1)n(\overline{n}n)^{r}\overline{n}$,\\
$\lbrack C]_{i}=(\overline{n}n)^{r}\overline{n}$,\\
$\lbrack C]_{i}=(n-1)(\overline{n}n)^{r}\overline{n}$,\\
$\lbrack C]_{i}=(n-1)(\overline{n}n)^{r}\overline{n}(\overline{n-1}).$%
\end{tabular}
\right.  \text{ if }i=n-1\text{ and }\left\{
\begin{tabular}
[c]{l}%
$\lbrack C]_{i}=n-1(n\overline{n})^{r}$,\\
$\lbrack C]_{i}=\overline{n}(n\overline{n})^{r}n$,\\
$\lbrack C]_{i}=(n-1)\overline{n}(n\overline{n})^{r}n$,\\
$\lbrack C]_{i}=(n\overline{n})^{r}n$,\\
$\lbrack C]_{i}=(n-1)(n\overline{n})^{r}n$,\\
$\lbrack C]_{i}=(n-1)(n\overline{n})^{r}n(\overline{n-1}).$%
\end{tabular}
\right.  \text{ if }i=n.
\]
where $(\overline{n}n)^{r}$ (resp. $(n\overline{n})^{r}$) is the word
containing the factor $\overline{n}n$ (resp.\ $n\overline{n}$) repeated $r$ times.
\end{proof}

Using Lemma \ref{Lem_dC_impl_Cadm} we derive immediately the

\begin{corollary}
A column $C$ of type $B$ or $D$ is admissible if and only if it can be split.
\end{corollary}

\begin{example}
\label{exa_splitC}From Example \ref{exam_splitting}, we obtain that $C$ is
admissible for $n=9$ and $C^{\prime}$ is admissible for $n=8.$
\end{example}

\begin{remark}
With the notations of the previous proposition, denote by $W_{n}/W_{\omega}$
the set of cosets of the Weyl group $W_{n}$ with respect to the stabilizer
$W_{\omega}$ of $\omega$ in $W_{n}$. Then we obtain a bijection $\tau$ between
the orbit \textsl{O}$_{\omega}$\textsl{\ }of $v_{\omega}$ in $B(\omega)$ under
the action of $W_{n}$ defined by (\ref{actionW}) and $W_{n}/W_{\omega}$. Using
Formulas (\ref{actionW}) it is easy to prove that \textsl{O}$_{\omega}$
consists of the vertices of $B(v_{\omega})$ without pair $(z,\overline{z})$.
Moreover if $C_{1},C_{2}$ are two columns such that $\mathrm{w(}C_{1}%
)=x_{1}\cdot\cdot\cdot x_{p},\mathrm{w(}C_{2})=y_{1}\cdot\cdot\cdot y_{p}%
\in\textsl{O}_{\omega}$, we have
\[
C_{1}\preceq C_{2}\Longleftrightarrow\tau_{\mathrm{w(}C_{1})}\vartriangleleft
_{\omega}\tau_{\mathrm{w(}C_{2})}%
\]
where $C_{1}\preceq C_{2}$ means that $x_{i}\preceq y_{i},i=1,...,p$ and
$"\vartriangleleft_{\omega}"$ denotes the projection of the Bruhat order on
$W_{n}/W_{\omega}.$ Then Proposition \ref{prop_imag_S2} may be regarded as a
version of Littelmann's labelling of $B(v_{\omega})$ by pairs $(\tau
_{\mathrm{w(}rC)},\tau_{\mathrm{w(}lC)})\in W_{n}/W_{\omega}\times
W_{n}/W_{\omega}$ satisfying $\tau_{\mathrm{w(}lC)}\vartriangleleft_{\omega
}\tau_{\mathrm{w(}rC)}$ \cite{Lit2}.
\end{remark}

\subsubsection{Orthogonal tableaux}

Every $\lambda\in\Omega_{+}^{B}$ has a unique decomposition of the form
$\lambda=\overset{n}{\underset{i=1}{\sum}}\lambda_{i}\omega_{i}^{B}%
$.\ Similarly, every $\lambda\in\Omega_{+}^{D}$ has a unique decomposition of
the form $\mathrm{(\ast)}$ $\lambda=\overset{n}{\underset{i=1}{\sum}}%
\lambda_{i}\omega_{i}^{D}$ or $\mathrm{(\ast\ast)}$ $\lambda=\lambda_{n}%
\bar{\omega}_{n}^{D}+\overset{n-1}{\underset{i=1}{\sum}}\lambda_{i}\omega
_{i}^{D} $ with $\lambda_{n}\neq0,$ where $(\lambda_{n},...,\lambda_{n}%
)\in\mathbb{N}^{n}.$ We will say that $(\lambda_{1},...,\lambda_{n})$ is the
positive decomposition of $\lambda\in\Omega_{+}$. Denote by $Y_{\lambda}$ the
Young diagram having $\lambda_{i}$ columns of height $i$ for $i=1,...,n$. If
$\lambda\in\Omega_{+}^{D},$ $Y_{\lambda}$ may not suffice to characterize the
weight $\lambda$ because a column diagram of length $n$ may be associated to
$\omega_{n}$ or to $\overline{\omega}_{n}$. In Subsection
\ref{sub_sec_Cor_in_Gn} we will need to attach to each dominant weight
$\lambda\in\Omega_{+}$ a combinatorial object $Y(\lambda)$. Moreover it will
be convenient to distinguish in $\mathrm{(\ast)}$ the cases where $\lambda
_{n}=0$ or $\lambda_{n}\neq0$. This leads us to set:
\begin{align}
\mathrm{(i)}  &  :Y(\lambda)=Y_{\lambda}\text{ if }\lambda\in\Omega_{+}%
^{B},\nonumber\\
\mathrm{(ii)}  &  :Y(\lambda)=(Y_{\lambda},+)\text{ in case }\mathrm{(\ast
)}\text{ with }\lambda_{n}\neq0,\nonumber\\
\mathrm{(iii)}  &  :Y(\lambda)=(Y_{\lambda},0)\text{ in case }\mathrm{(\ast
)}\text{ with }\lambda_{n}=0,\label{Def_Y(lambda)}\\
\mathrm{(iv)}  &  :Y(\lambda)=(Y_{\lambda},-)\text{ in case }\mathrm{(\ast
\ast).}\nonumber
\end{align}
When $\lambda\in\Omega_{+}^{D},$ $Y(\lambda)$ may be regarded as the
generalization of the notion of the shape of type $A$ associated to a dominant
weight. Now write
\begin{align*}
v_{\lambda}^{B}  &  =(v_{\omega_{1}^{B}})^{\otimes\lambda_{1}}\otimes
\cdot\cdot\cdot\otimes(v_{\omega_{n}^{B}})^{\otimes\lambda_{n}}\text{ in case
}\mathrm{(i),}\\
v_{\lambda}^{D}  &  =(v_{\omega_{1}^{D}})^{\otimes\lambda_{1}}\otimes
\cdot\cdot\cdot\otimes(v_{\omega_{n}^{D}})^{\otimes\lambda_{n}}\text{ in
case$\mathrm{(ii)}$,}\\
v_{\lambda}^{D}  &  =(v_{\omega_{1}^{D}})^{\otimes\lambda_{1}}\otimes
\cdot\cdot\cdot\otimes(v_{\omega_{n-1}^{D}})^{\otimes\lambda_{n-1}}\text{ in
case $\mathrm{(iii)}$ and}\\
v_{\lambda}^{D}  &  =(v_{\omega_{1}^{D}})^{\otimes\lambda_{1}}\otimes
\cdot\cdot\cdot\otimes(v_{\overline{\omega}_{n}^{D}})^{\otimes\lambda_{n}%
}\text{ in case $\mathrm{(iv).}$}%
\end{align*}
Then $v_{\lambda}^{B}$ and $v_{\lambda}^{D}$ are highest weight vertices of
$G_{n}^{B}$ and $G_{n}^{D}.$ Moreover $B(v_{\lambda}^{B})$ and $B(v_{\lambda
}^{D})$ are isomorphic to $B^{B}(\lambda)$ and $B^{D}(\lambda)$.

A tabloid $\tau$ of type $B$ (resp. $D$) is a Young diagram whose columns are
filled to give columns of type $B$ (resp. $D$). If $\tau=C_{1}\cdot\cdot\cdot
C_{r}$, we write \textrm{w}$(T)=$\textrm{w}$(C_{r})\cdot\cdot\cdot$%
\textrm{w}$(C_{1})$ for the reading of $\tau$.

\begin{definition}
\label{defKN} \ \ \ \ \ \ 

\begin{itemize}
\item  Consider $\lambda\in\Omega_{+}^{B}$. A tabloid $T$ of type $B$ is an
orthogonal tableau of shape $Y(\lambda)$ and type $B$ if $\mathrm{w(}%
T\mathrm{)}\in$ $B(v_{\lambda}^{B})$.

\item  Consider $\lambda\in\Omega_{+}^{D}$. A tabloid $T$ of type $D$ is an
orthogonal tableau of shape $Y(\lambda)$ and type $D$ if $\mathrm{w(}%
T\mathrm{)}\in$ $B(v_{\lambda}^{D})$.
\end{itemize}
\end{definition}

The orthogonal tableaux of a given shape form a single connected component of
$G_{n}$, hence two orthogonal tableaux whose readings occur at the same place
in two isomorphic connected components of $G_{n}$ are equal. The shape of an
orthogonal tableau $T$ of type $D$ may be regarded as a pair $[O_{T}%
,\varepsilon_{T}]$ where $O_{T}$ is a Young diagram and $\varepsilon_{T}%
\in\{-,0,+\}$. The $\{-,0,+\}$ part of this shape can be read off directly on
$T$.\ Indeed $\varepsilon=0$ if $T$ does not contain a column of height $n.$
Otherwise write $\mathrm{w(}C_{1})=x_{1}\cdot\cdot\cdot x_{n}$ for the reading
of the first column of $T.$ Since it is admissible, $C_{1}$ contains at least
a letter, say $x_{k}$ of $\{n,\overline{n}\}.$ Then $\varepsilon$ is given by
the parity of $n-k$ according to Proposition \ref{prop_KNf}.

Consider $\tau=C_{1}C_{2}\cdot\cdot\cdot C_{r}$ a tabloid whose columns are
admissible. The split form of $\tau$ is the tabloid obtained by splitting each
column of $\tau$.\ We write $\mathrm{spl}(\tau)=(lC_{1}rC_{1})(lC_{2}%
rC_{2})\cdot\cdot\cdot(lC_{r}rC_{r})$.\ With the notations of Proposition
\ref{prop_imag_S2}, we will have $\mathrm{w}(\mathrm{spl}(T))=S_{2}%
\mathrm{w(}C_{r})\cdot\cdot\cdot S_{2}\mathrm{w(}C_{1}) $%
.\ Kashiwara-Nakashima's combinatorial description \cite{KN} of an orthogonal
tableau $T$ is based on the enumeration of configurations that should not
occur in two adjacent columns of $T$. Considering its split form
$\mathrm{spl}(T)$, this description becomes more simple because the columns of
$\mathrm{spl}(T)$ does not contain any pair $(z,\overline{z})$.

\begin{lemma}
\label{lem_split_tab}Let $T=C_{1}C_{2}\cdot\cdot\cdot C_{r}$ be a tabloid
whose columns are admissible. Then $T$ is an orthogonal tableau if and only if
$\mathrm{spl}(T)$ is an orthogonal tableau.
\end{lemma}

\begin{proof}
Suppose first that $\mathrm{w(}T)$ is a highest weight vertex of weight
$\lambda$. Then, by Corollary \ref{cor_strech}, $\mathrm{w(spl}(T))$ is a
highest weight vertex of weight $2\lambda$. If $T$ is an orthogonal tableau,
$\mathrm{w(}T)=v_{\lambda}$ and we have $\mathrm{w(spl}(T))=v_{2\lambda}$. So
$\mathrm{spl}(T)$ is an orthogonal tableau. Conversely, if $\mathrm{spl}(T)$
is an orthogonal tableau, $\mathrm{w(spl}(T))=S_{2}\mathrm{w(}C_{r})\cdot
\cdot\cdot S_{2}\mathrm{w(}C_{1})$ is a highest weight vertex of weight
$2\lambda$ by Corollary \ref{cor_strech}. Hence we have $\mathrm{w(spl}%
(T))=v_{2\lambda}$ because there exists only one orthogonal tableau of highest
weight $2\lambda$. So $\mathrm{w(}T)=v_{\lambda}$. In the general case, denote
by $T_{0}$ the tableau such that $\mathrm{w(}T_{0})$ is the highest weight
vertex of the connected component of $G_{n}$ containing $\mathrm{w(}T)$. Then
$\mathrm{w(spl}(T_{0}))$ is the highest weight vertex of the connected
component containing $\mathrm{w(spl}(T))$ and the following assertions are
equivalent:
\begin{align*}
&  \text{$\mathrm{(i)}$\textrm{\ }}\mathrm{spl}(T)\text{ is an orthogonal
tableau,}\\
&  \text{$\mathrm{(ii)}$\textrm{\ }}\mathrm{spl}(T_{0})\text{ is orthogonal
tableau,}\\
&  \text{$\mathrm{(iii)}$\textrm{\ }}T_{0}\text{ is orthogonal tableau,}\\
&  \text{$\mathrm{(iv)}$\textrm{\ }}T\text{ is orthogonal tableau.}%
\end{align*}
\end{proof}

\begin{definition}
Let $\tau=C_{1}C_{2}$ be a tabloid with two admissible columns $C_{1}$ and
$C_{2}$. We set:

\begin{itemize}
\item $C_{1}\preceq C_{2}$ when $h(C_{1})\geq h(C_{2})$ and the rows of
$C_{1}C_{2}$ are weakly increasing from left to right,

\item $C_{1}\trianglelefteq C_{2}$ when $rC_{1}\preceq lC_{2}.$
\end{itemize}
\end{definition}

\begin{definition}
\label{Def_b_conf}(Kashiwara-Nakashima)

\noindent Let $C_{1}=%
\begin{tabular}
[c]{|l|}\hline
$x_{1}$\\\hline
$\cdot$\\\hline
$\cdot$\\\hline
$x_{N}$\\\hline
\end{tabular}
$ and $C_{2}=%
\begin{tabular}
[c]{|l|}\hline
$y_{1}$\\\hline
$\cdot$\\\hline
$\cdot$\\\hline
$y_{N}$\\\hline
\end{tabular}
$ be admissible columns of type $D$ and $p,q,r,s$ integers satisfying $1\leq
p\leq q<r\leq s\leq M$.

$C_{1}C_{2}$ contains an $a$-odd-configuration (with $a\notin\{\overline
{n},n\}$) when:

\begin{itemize}
\item $a=x_{p},\overline{n}=x_{r}$ are letters of $C_{1}$ and $\overline
{a}=y_{s},n=y_{q}$ letters of $C_{2}$ such that $r-q+1$ is odd
\end{itemize}

or

\begin{itemize}
\item $a=x_{p},n=x_{r}$ are letters of $C_{1}$ and $\overline{a}%
=y_{s},\overline{n}=y_{q}$ letters of $C_{2}$ such that $r-q+1$ is odd
\end{itemize}

$C_{1}C_{2}$ contains an $a$-even-configuration (with $a\notin\{\overline
{n},n\}$) when:

\begin{itemize}
\item $a=x_{p},n=x_{r}$ are letters of $C_{1}$ and $\overline{a}=y_{s}%
,n=y_{q}$ letters of $C_{2}$ such that $r-q+1$ is even
\end{itemize}

or

\begin{itemize}
\item $a=x_{p},\overline{n}=x_{r}$ are letters of $C_{1}$ and $\overline
{a}=y_{s},\overline{n}=y_{q}$ letters of $C_{2}$ such that $r-q+1$ is even
\end{itemize}

Then we denote by $\mu(a)$ the positive integer defined by:
\[
\mu(a)=s-p
\]
\end{definition}

\begin{theorem}
\ \ \ \ \ \label{TH_KN}

\noindent$\mathrm{(i)}$ Consider $C_{1},C_{2},...,C_{r}$ some admissible
columns of type $B$. Then the tabloid $T=C_{1}C_{2}\cdot\cdot\cdot C_{r}$ is
an orthogonal tableau if and only if $C_{i}\trianglelefteq C_{i+1}$ for $i=1,...,r-1.$

\noindent$\mathrm{(ii)}$ Consider $C_{1},C_{2},...,C_{r}$ some admissible
columns of type $D$. Then the tabloid $T=C_{1}C_{2}\cdot\cdot\cdot C_{r}$ is
an orthogonal tableau if and only if, $C_{i}\trianglelefteq C_{i+1}$ for
$i=1,...,r-1$, and $rC_{i}lC_{i+1}$ does not contain an $a$-configuration
(even or odd) such that $\mu(a)=n-a$.
\end{theorem}

\begin{proof}
Kashiwara and Nakashima describe an orthogonal tableau $T$ by listing the
configurations that should not occur in two adjacent columns of $T$. If we
except the $a$-configurations even or odd, these configurations disappear in
$\mathrm{spl}(T)$ because $\mathrm{spl}(T)$ does not contain a column with a
pair $(z,\overline{z}).$ Hence the theorem follows from Lemma
\ref{lem_split_tab} and Theorems 5.7.1, 6.7.1 of \cite{KN}.
\end{proof}

\begin{example}
Suppose $n=4.$ Then

\noindent$T=%
\begin{tabular}
[c]{|l|ll}\hline
$\mathtt{3}$ & $\mathtt{3}$ & \multicolumn{1}{|l|}{$\mathtt{4}$}\\\hline
$\mathtt{4}$ & $\mathtt{0}$ & \multicolumn{1}{|l|}{$\mathtt{\bar{4}}$}\\\hline
$\mathtt{0}$ & $\mathtt{\bar{2}}$ & \multicolumn{1}{|l}{}\\\cline{1-2}%
$\mathtt{0}$ &  & \\\cline{1-1}%
\end{tabular}
$ is an orthogonal tableau of type $B$ because $\mathrm{spl}(T)=%
\begin{tabular}
[c]{|l|l|llll}\hline
$\mathtt{1}$ & $\mathtt{3}$ & $\mathtt{3}$ & \multicolumn{1}{|l}{$\mathtt{3}$}%
& \multicolumn{1}{|l}{$\mathtt{3}$} & \multicolumn{1}{|l|}{$\mathtt{4}$%
}\\\hline
$\mathtt{2}$ & $\mathtt{4}$ & $\mathtt{4}$ & \multicolumn{1}{|l}{$\mathtt{\bar
{4}}$} & \multicolumn{1}{|l}{$\mathtt{\bar{4}}$} &
\multicolumn{1}{|l|}{$\mathtt{\bar{3}}$}\\\hline
$\mathtt{3}$ & $\mathtt{\bar{2}}$ & $\mathtt{\bar{2}}$ &
\multicolumn{1}{|l}{$\mathtt{\bar{2}}$} & \multicolumn{1}{|l}{} &
\\\cline{1-2}\cline{1-4}%
$\mathtt{4}$ & $\mathtt{\bar{1}}$ &  &  &  & \\\cline{1-2}%
\end{tabular}
\vspace{0.2cm}.$ But
\begin{tabular}
[c]{|l|l|}\hline
$\mathtt{3}$ & $\mathtt{\bar{4}}$\\\hline
$\mathtt{\bar{4}}$ & $\mathtt{\bar{3}}$\\\hline
\end{tabular}
is not orthogonal of type $D$ because it contains a $3$-even configuration
with $\mu(3)=1$.
\end{example}

\subsection{\label{subsec_monoids}Plactic monoids for types $B_{n}$ and
$D_{n}$}

\begin{definition}
\label{def_sam_plac}Let $w_{1}$ and $w_{2}$ be two words on $\mathcal{B}_{n}$
(resp. $\mathcal{D}_{n}$) . We write $w_{1}\overset{B}{\sim}w_{2}$ (resp.
$w_{1}\overset{D}{\sim}w_{2}$) when these two words occur at the same place in
two isomorphic connected components of the crystal $G_{n}^{B}$ (resp.
$G_{n}^{D}$).
\end{definition}

The definition of the orthogonal tableaux implies that for any word
$w\in\mathcal{B}_{n}^{\ast}$ (resp. $w\in\mathcal{D}_{n}^{\ast}$) there exists
a unique orthogonal tableau $P^{B}(w)$ (resp. $P^{D}(w)$) such that
$w\sim\mathrm{w(}P(w))$. So the sets $\mathcal{B}_{n}^{\ast}/\overset{B}{\sim
}$ and $\mathcal{D}_{n}^{\ast}/\overset{D}{\sim}$ can be identified
respectively with the sets of orthogonal tableaux of type $B$ and $D$. Our aim
is now to show that $\overset{B}{\sim}$ and $\overset{D}{\sim}$ are in fact
congruencies $\overset{B}{\equiv}$ and $\overset{D}{\equiv}$ so that
$\mathcal{B}_{n}^{\ast}/\overset{B}{\sim}$ and $\mathcal{D}_{n}^{\ast
}/\overset{D}{\sim}$ are in a natural way endowed with a multiplication.

\begin{definition}
The monoid $Pl(B_{n})$ is the quotient of the free monoid $\mathcal{B}%
_{n}^{\ast}$ by the relations:

$R_{1}^{B}:$ If $x\neq\overline{z}$ and $x\prec y\prec z:$%
\[
yzx\overset{B}{\equiv}yxz\text{ \ and \ }xzy\overset{B}{\equiv}zxy\text{.}%
\]

$R_{2}^{B}:$ If $x\neq\overline{y}$ and $x\prec y:$%
\[
xyx\overset{B}{\equiv}xxy\text{ for }x\neq0\text{ \ and \ }xyy\overset
{B}{\equiv}yxy\text{ for }y\neq0.
\]

$R_{3}^{B}:$If $1\prec x\preceq n$ and $x\preceq y\preceq\overline{x}:$%

\[
y(\overline{x-1})(x-1)\overset{B}{\equiv}yx\overline{x}\text{, \ and
\ }x\overline{x}y\overset{B}{\equiv}(\overline{x-1})(x-1)y,
\]%
\[
0\overline{n}n\equiv\overline{n}n0.
\]

$R_{4}^{B}:$If $x\preceq n:$%
\[
00x\overset{B}{\equiv}0x0\text{ \ and \ }0\overline{x}0\overset{B}{\equiv
}\overline{x}00.
\]

$R_{5}^{B}:$ Let $w=\mathrm{w}(C)$ be a non admissible column word each strict
factor of which is admissible. When $C$ satisfies the assertion $\mathrm{(i)}$
of Remark \ref{not_N(z)}, let $z$ be the lowest unbarred letter of $w$ such
that the pair $(z,\overline{z})$ occurs in $w$ and $N(z)>z$, otherwise set
$z=0$. Then $w\overset{B}{\equiv}\widetilde{w}$ where $\widetilde{w}$ is the
column word obtained by erasing the pair $(z,\overline{z})$ in $w$ if
$z\preceq n,$ by erasing $0$ otherwise.
\end{definition}

\begin{definition}
The monoid $Pl(D_{n})$ is the quotient of the free monoid $\mathcal{D}%
_{n}^{\ast}$ by the relations:

$R_{1}:$ If $x\neq\overline{z}$
\[
yzx\overset{D}{\equiv}yxz\text{ for }x\preceq y\prec z\text{ \ and
\ }xzy\overset{D}{\equiv}zxy\text{ for }x\prec y\preceq z.
\]

$R_{2}:$ If $1\prec x\prec n$ and $x\preceq y\preceq\overline{x}$%
\[
y(\overline{x-1})(x-1)\overset{D}{\equiv}yx\overline{x}\text{ \ and
\ }x\overline{x}y\overset{D}{\equiv}(\overline{x-1})(x-1)y.
\]

$R_{3}^{D}:$ If $x\preceq n-1:$%
\[
\left\{
\begin{tabular}
[c]{l}%
$\overline{n}\,\overline{x}n\overset{D}{\equiv}\overline{x}\,\overline{n}n$\\
$n\,\overline{x}\,\overline{n}\overset{D}{\equiv}\overline{x}\,n\overline{n}$%
\end{tabular}
\right.  \text{ and }\left\{
\begin{tabular}
[c]{l}%
$\overline{n}nx\overset{D}{\equiv}\overline{n}xn$\\
$n\overline{n}x\overset{D}{\equiv}nx\overline{n}$%
\end{tabular}
\right.  .
\]

$R_{4}^{D}:$%
\[
\left\{
\begin{tabular}
[c]{l}%
$n\overline{n}\,\overline{n}\overset{D}{\equiv}\overline{(n-1)}%
\,(n-1)\overline{n}$\\
$\overline{n}\,nn\overset{D}{\equiv}\overline{(n-1)}\,(n-1)n$%
\end{tabular}
\right.  \text{ and }\left\{
\begin{tabular}
[c]{l}%
$\overline{n}(\overline{n-1})(n-1)\overset{D}{\equiv}\overline{n}%
\,\overline{n}n$\\
$n(\overline{n-1})(n-1)\overset{D}{\equiv}nn\,\overline{n}$%
\end{tabular}
\right.  .
\]

$R_{5}^{D}:$ Consider $w$ a non admissible column word each strict factor of
which is admissible. Let $z$ be the lowest unbarred letter such that the pair
$(z,\overline{z})$ occurs in $w$ and $N(z)>z$ (see Remark \ref{not_N(z)}).
Then $w\overset{D}{\equiv}\widetilde{w}$ where $\widetilde{w}$ is the column
word obtained by erasing the pair $(z,\overline{z})$ in $w$ if $z\prec n$, by
erasing a pair $(n,\overline{n})$ of consecutive letters otherwise.
\end{definition}

The relations $R_{5}^{B}$ and $R_{5}^{D}$ are called the contraction
relations.\ When the letter $0$ or a pair $(n,\overline{n})$ disappears, we
have $l(C)=n+1$ and in $R_{5}^{D}$ the word $\widetilde{w}$ does not depend on
the factor $n\overline{n}$ or $\overline{n}n$ erased. Moreover $\widetilde{w}$
is an admissible column word. Note that $w_{1}\equiv w_{2}$ implies
$d(w_{1})=d(w_{2})$, that is, $\equiv$ is compatible with the grading given by
$d$.

\begin{theorem}
\label{th_Psymbol}Given two words $w_{1}$ and $w_{2}$%
\begin{equation}
w_{1}\sim w_{2}\Longleftrightarrow w_{1}\equiv w_{2}\Longleftrightarrow
P(w_{1})=P(w_{2}) \label{good_rela}%
\end{equation}
\end{theorem}

\noindent This theorem is proved in the same way as in the symplectic case
\cite{Lec}, and we will only sketch the arguments. Note first that we have%
\[
w_{1}\sim w_{2}\Longleftrightarrow P(w_{1})=P(w_{2})
\]
immediately from the definition of $P.$ For any word $w$ occurring in the left
hand side of a relation $R_{1}^{B},...,R_{4}^{B}$ (resp.\ $R_{1}^{D}%
,...,R_{4}^{D}$), write $\xi^{B}(w)$ (resp. $\xi^{D}(w)$) for the word
occurring in the right hand side of this relation. Similarly for $p=1,...,n$
and $w$ a word of length $p+1$ occurring in the left hand side of $R_{5}^{B}$
(resp. $R_{5}^{D}$), denote by $\xi_{p}^{B}(w)$ (resp. $\xi_{p}^{D}(w)$) the
word occurring in the right hand side of this relation. By using similar
arguments to those of \cite{Lec}, we obtain the followings assertions:

\begin{itemize}
\item  The map $\xi^{B}:w\longmapsto\xi(w)$ is the crystal isomorphism from
$B^{B}(121)$ to $B^{B}(112)$.

\item  If $n>2,$ the map $\xi^{D}:w\longmapsto\xi(w)$ is the crystal
isomorphism from $B^{D}(121)$ to $B^{D}(112)$ otherwise $\xi^{D}$ is the
crystal isomorphism from $B^{D}(121)\cup B^{D}(1\bar{2}1)$ to $B^{D}(112)\cup
B^{D}(11\bar{2})$.

\item  For $p=2,...,n-1$, $\xi_{p}:w\longmapsto\xi_{p}(w)$ is the crystal
isomorphism from $B(12\cdot\cdot\cdot p\overline{p})$ to $B(12\cdot\cdot\cdot
p-1)$.

\item  The map $\xi_{n}^{B}:w\longmapsto\xi_{n}^{B}(w)$ is the crystal
isomorphism from $B^{B}(12\cdot\cdot\cdot n\overline{n})\cup B^{B}%
(12\cdot\cdot\cdot n0)$ to $B^{B}(12\cdot\cdot\cdot n-1)\cup B^{B}%
(12\cdot\cdot\cdot n)$.

\item  The words $w$ of length $n+1$ occurring in the left hand side of
$R_{5}^{D}$ are the vertices of $B^{D}(12\cdot\cdot\cdot n\overline{n})\cup
B^{D}(12\cdot\cdot\cdot\overline{n}n)$. Moreover the restriction of the map
$\xi_{n}^{D}:w\longmapsto\xi_{n}^{D}(w)$ to $B^{D}(12\cdot\cdot\cdot
n\overline{n})$ (resp. to $B^{D}(12\cdot\cdot\cdot\overline{n}n)$) is the
crystal isomorphism from $B^{D}(12\cdot\cdot\cdot n\overline{n})$ (resp.
$B^{D}(12\cdot\cdot\cdot\overline{n}n)$) to $B^{D}(12\cdot\cdot\cdot n-1)$.
\end{itemize}%

\begin{gather*}%
\begin{tabular}
[c]{lllllllll}%
&  &  &  & $121$ &  &  &  & \\
&  &  & $\overset{\text{1}}{\swarrow}$ &  & $\overset{\text{2}}{\searrow}$ &
&  & \\
&  & $122$ &  &  &  & $101$ &  & \\
& $\overset{\text{2}}{\swarrow}$ &  &  &  &  & {\tiny 1}$\downarrow$ &
$\overset{\text{2}}{\searrow}$ & \\
$102$ &  &  &  &  &  & $201$ &  & $1\bar{2}1$\\
{\tiny 2}$\downarrow$ &  &  &  &  & $\overset{\text{1}}{\swarrow}$ &
{\tiny 2}$\downarrow$ &  & {\tiny 1}$\downarrow$\\
$1\bar{2}2$ &  &  &  & $202$ &  & $001$ &  & $2\bar{2}1$\\
{\tiny 2}$\downarrow$ & $\overset{\text{1}}{\searrow}$ &  &  & {\tiny 2}%
$\downarrow$ & $\overset{\text{1}}{\swarrow}$ & {\tiny 2}$\downarrow$ &  &
{\tiny 1}$\downarrow$\\
$1\bar{2}0$ &  & $2\bar{2}2$ &  & $002$ &  & $0\bar{2}1$ &  & $2\bar{1}1$\\
{\tiny 2}$\downarrow$ & $\overset{\text{1}}{\searrow}$ & {\tiny 2}$\downarrow$%
&  & {\tiny 2}$\downarrow$ &  & {\tiny 1}$\downarrow$ & $\overset{\text{2}%
}{\swarrow}$ & {\tiny 1}$\downarrow$\\
$1\bar{2}\bar{2}$ &  & $2\bar{2}0$ &  & $0\bar{2}2$ &  & $0\bar{1}1$ &  &
$2\bar{1}2$\\
& $\overset{\text{2}}{\swarrow}$ & {\tiny 1}$\downarrow$ &  & {\tiny 2}%
$\downarrow$ &  & {\tiny 2}$\downarrow$ & $\overset{\text{1}}{\searrow}$ &
{\tiny 2}$\downarrow$\\
$2\bar{2}\bar{2}$ &  & $2\bar{1}0$ &  & $0\bar{2}0$ &  & $\bar{2}\bar{1}1$ &
& $0\bar{1}2$\\
{\tiny 1}$\downarrow$ &  & {\tiny 2}$\downarrow$ & $\overset{\text{1}%
}{\swarrow}$ & {\tiny 2}$\downarrow$ &  &  & $\overset{\text{1}}{\searrow}$ &
{\tiny 2}$\downarrow$\\
$2\bar{1}\bar{2}$ &  & $0\bar{1}0$ &  & $0\bar{2}\bar{2}$ &  &  &  & $\bar
{2}\bar{1}2$\\
{\tiny 1}$\downarrow$ &  & {\tiny 2}$\downarrow$ & $\overset{\text{1}%
}{\swarrow}$ &  &  &  &  & {\tiny 2}$\downarrow$\\
$2\bar{1}\bar{1}$ &  & $0\bar{1}\bar{2}$ &  &  &  &  &  & $\bar{2}\bar{1}0$\\
& $\overset{\text{2}}{\searrow}$ & {\tiny 1}$\downarrow$ &  &  &  &  &
$\overset{\text{2}}{\swarrow}$ & \\
&  & $0\bar{1}\bar{1}$ &  &  &  & $\bar{2}\bar{1}\bar{2}$ &  & \\
&  &  & $\overset{\text{2}}{\searrow}$ &  & $\overset{\text{1}}{\swarrow}$ &
&  & \\
&  &  &  & $\bar{2}\bar{1}\bar{1}$ &  &  &  & \\
&  &  &  &  &  &  &  &
\end{tabular}
\text{ \ \ \ \ \ }
\begin{tabular}
[c]{lllllllll}%
&  &  &  & $112$ &  &  &  & \\
&  &  & $\overset{\text{1}}{\swarrow}$ &  & $\overset{\text{2}}{\searrow}$ &
&  & \\
&  & $212$ &  &  &  & $110$ &  & \\
& $\overset{\text{2}}{\swarrow}$ &  &  &  &  & {\tiny 1}$\downarrow$ &
$\overset{\text{2}}{\searrow}$ & \\
$012$ &  &  &  &  &  & $210$ &  & $11\bar{2}$\\
{\tiny 2}$\downarrow$ &  &  &  &  & $\overset{\text{1}}{\swarrow}$ &
{\tiny 2}$\downarrow$ &  & {\tiny 1}$\downarrow$\\
$\bar{2}12$ &  &  &  & $220$ &  & $010$ &  & $21\bar{2}$\\
{\tiny 2}$\downarrow$ & $\overset{\text{1}}{\searrow}$ &  &  & {\tiny 2}%
$\downarrow$ & $\overset{\text{1}}{\swarrow}$ & {\tiny 2}$\downarrow$ &  &
{\tiny 1}$\downarrow$\\
$\bar{2}10$ &  & $\bar{1}12$ &  & $020$ &  & $01\bar{2}$ &  & $22\bar{2}$\\
{\tiny 2}$\downarrow$ & $\overset{\text{1}}{\searrow}$ & {\tiny 2}$\downarrow$%
&  & {\tiny 2}$\downarrow$ &  & {\tiny 1}$\downarrow$ & $\overset{\text{2}%
}{\swarrow}$ & {\tiny 1}$\downarrow$\\
$\bar{2}1\bar{2}$ &  & $\bar{1}10$ &  & $\bar{2}20$ &  & $02\bar{2}$ &  &
$22\bar{1}$\\
& $\overset{\text{2}}{\swarrow}$ & {\tiny 1}$\downarrow$ &  & {\tiny 2}%
$\downarrow$ &  & {\tiny 2}$\downarrow$ & $\overset{\text{1}}{\searrow}$ &
{\tiny 2}$\downarrow$\\
$\bar{1}1\bar{2}$ &  & $\bar{1}20$ &  & $\bar{2}00$ &  & $\bar{2}2\bar{2}$ &
& $02\bar{1}$\\
{\tiny 1}$\downarrow$ &  & {\tiny 2}$\downarrow$ & $\overset{\text{1}%
}{\swarrow}$ & {\tiny 2}$\downarrow$ &  &  & $\overset{\text{1}}{\searrow}$ &
{\tiny 2}$\downarrow$\\
$\bar{1}2\bar{2}$ &  & $\bar{1}00$ &  & $\bar{2}0\bar{2}$ &  &  &  & $\bar
{2}2\bar{1}$\\
{\tiny 1}$\downarrow$ &  & {\tiny 2}$\downarrow$ & $\overset{\text{1}%
}{\swarrow}$ &  &  &  &  & {\tiny 2}$\downarrow$\\
$\bar{1}2\bar{1}$ &  & $\bar{1}0\bar{2}$ &  &  &  &  &  & $\bar{2}0\bar{1}$\\
& $\overset{\text{2}}{\searrow}$ & {\tiny 1}$\downarrow$ &  &  &  &  &
$\overset{\text{2}}{\swarrow}$ & \\
&  & $\bar{1}0\bar{1}$ &  &  &  & $\bar{2}\bar{2}\bar{1}$ &  & \\
&  &  & $\overset{\text{2}}{\searrow}$ &  & $\overset{\text{1}}{\swarrow}$ &
&  & \\
&  &  &  & $\bar{1}\bar{2}\bar{1}$ &  &  &  & \\
&  &  &  &  &  &  &  &
\end{tabular}
\\
\text{The crystals }B^{B}(121)\text{ and }B^{B}(112)\text{ in }G_{2}^{B}%
\end{gather*}

\bigskip%

\begin{gather*}%
\begin{tabular}
[c]{lll}%
$121$ &  & \\
{\tiny 1}$\downarrow$ & $\overset{\text{2}}{\searrow}$ & \\
$122$ &  & $\bar{2}21$\\
{\tiny 2}$\downarrow$ & $\overset{\text{1}}{\swarrow}$ & {\tiny 2}$\downarrow
$\\
$\bar{2}22$ &  & $\bar{2}\bar{1}1$\\
{\tiny 2}$\downarrow$ & $\overset{\text{1}}{\swarrow}$ & {\tiny 2}$\downarrow
$\\
$\bar{2}\bar{1}2$ &  & $\bar{2}\bar{1}\bar{2}$\\
{\tiny 2}$\downarrow$ & $\overset{\text{1}}{\swarrow}$ & \\
$\bar{2}\bar{1}\bar{1}$ &  &
\end{tabular}
\text{ \ \ \ \ \ \ \ }
\begin{tabular}
[c]{lll}%
$112$ &  & \\
{\tiny 1}$\downarrow$ & $\overset{\text{2}}{\searrow}$ & \\
$212$ &  & $\bar{2}12$\\
{\tiny 2}$\downarrow$ & $\overset{\text{1}}{\swarrow}$ & {\tiny 2}$\downarrow
$\\
$\bar{1}12$ &  & $\bar{2}\bar{2}2$\\
{\tiny 2}$\downarrow$ & $\overset{\text{1}}{\swarrow}$ & {\tiny 2}$\downarrow
$\\
$\bar{1}\bar{2}2$ &  & $\bar{2}\bar{2}\bar{1}$\\
{\tiny 2}$\downarrow$ & $\overset{\text{1}}{\swarrow}$ & \\
$\bar{1}\bar{2}\bar{1}$ &  &
\end{tabular}
\\
\text{The crystals }B^{D}(121)\text{ and }B^{D}(112)\text{ in }G_{2}^{D}%
\end{gather*}%

\begin{gather*}%
\begin{tabular}
[c]{lll}%
$1\bar{2}1$ &  & \\
{\tiny 1}$\downarrow$ & $\overset{\text{2}}{\searrow}$ & \\
$2\bar{2}1$ &  & $1\bar{2}\bar{2}$\\
{\tiny 1}$\downarrow$ & $\overset{\text{2}}{\searrow}$ & {\tiny 1}$\downarrow
$\\
$2\bar{1}1$ &  & $2\bar{2}\bar{2}$\\
{\tiny 1}$\downarrow$ &  & {\tiny 1}$\downarrow$\\
$2\bar{1}2$ &  & $2\bar{1}\bar{2}$\\
{\tiny 2}$\downarrow$ & $\overset{\text{1}}{\swarrow}$ & \\
$2\bar{1}\bar{1}$ &  &
\end{tabular}
\text{ \ \ \ \ \ \ }
\begin{tabular}
[c]{lll}%
$11\bar{2}$ &  & \\
{\tiny 1}$\downarrow$ & $\overset{\text{2}}{\searrow}$ & \\
$21\bar{2}$ &  & $\bar{2}1\bar{2}$\\
{\tiny 1}$\downarrow$ & $\overset{\text{2}}{\searrow}$ & {\tiny 1}$\downarrow
$\\
$22\bar{2}$ &  & $\bar{1}1\bar{2}$\\
{\tiny 1}$\downarrow$ &  & {\tiny 1}$\downarrow$\\
$22\bar{1}$ &  & $\bar{1}2\bar{2}$\\
{\tiny 2}$\downarrow$ & $\overset{\text{1}}{\swarrow}$ & \\
$\bar{1}2\bar{1}$ &  &
\end{tabular}
\\
\text{The crystals }B^{D}(1\bar{2}1)\text{ and }B^{D}(11\bar{2})\text{ in
}G_{2}^{D}%
\end{gather*}

By (\ref{TENS1}) and (\ref{TENS2}), this implies that the plactic relations
above are compatible with Kashiwara's operators, that is, for any words
$w_{1}$ and $w_{2}$ such that $w_{1}\equiv w_{2}$ one has:
\begin{equation}
\left\{
\begin{tabular}
[c]{l}%
$\widetilde{e}_{i}(w_{1})\equiv\widetilde{e}_{i}(w_{2})\text{ and }%
\varepsilon_{i}(w_{1})=\varepsilon_{i}(w_{2})$\\
$\widetilde{f}_{i}(w_{1})\equiv\widetilde{f}_{i}(w_{2})\text{ and }\varphi
_{i}(w_{1})=\varphi_{i}(w_{2}).$%
\end{tabular}
\right.  \label{fonda_compatib}%
\end{equation}
Hence:%

\[
w_{1}\equiv w_{2}\Longrightarrow w_{1}\sim w_{2}.
\]
To obtain the converse we show that for any highest weight vertex $w^{0}$%
\begin{equation}
\mathrm{w(}P(w^{0}))\equiv w^{0}. \label{cong_on_HWV}%
\end{equation}
This follows by induction on $\mathrm{l}(w^{0})$.\ When $\mathrm{l}(w^{0})=1$,
$\mathrm{w(}P(w^{0}))=w^{0}$. By writing $w^{0}=v^{0}x^{0}$, it is possible
(see the proof of Lemma 3.2.6 in \cite{Lec}) to show that $\mathrm{w(}%
P(w^{0}))$ may be obtained from the word $\mathrm{w(}P(v^{0}))x^{0}$ by
applying only Knuth relations and contractions relations of type $12\cdot
\cdot\cdot r\overline{p}\equiv12\cdot\cdot\cdot\widehat{p}\cdot\cdot\cdot r$
with $p\leq r\leq n$ (the hat means removal the letter $p$).

From (\ref{cong_on_HWV}), we obtain that two highest weight vertices
$w_{1}^{0}$ and $w_{2}^{0}$ with the same weight $\lambda$ verify $w_{1}%
^{0}\equiv w_{2}^{0}$. Indeed there is only one orthogonal tableau whose
reading is a highest vertex of weight $\lambda$. Now suppose that $w_{1}\sim
w_{2}$ and denote by $w_{1}^{0}$ and $w_{2}^{0}$ the highest weight vertices
of $B(w_{1})$ and $B(w_{2})$. We have $w_{1}^{0}\equiv w_{2}^{0}$. Set
$w_{1}=\widetilde{F}\,w_{1}^{0}$ where $\widetilde{F}$ is a product of
Kashiwara's operators $\widetilde{f}_{i}$, $i=1,...,n$. Then $w_{2}%
=\widetilde{F}\,w_{2}^{0}$ because $w_{1}\sim w_{2}$. So by
(\ref{fonda_compatib}) we obtain
\[
w_{1}^{0}\equiv w_{2}^{0}\Longrightarrow\widetilde{F}\,w_{1}^{0}%
\equiv\widetilde{F}\,w_{2}^{0}\Longrightarrow w_{1}\equiv w_{2}.
\]

\subsection{A bumping algorithm for types $B$ and $D$}

Now we are going to see how the orthogonal tableau $P(w)$ may be computed for
each vertex $w$ by using an insertion scheme analogous to bumping algorithm
for type $A$. As a first step, we describe $P(w)$ when $w=\mathrm{w}(C)x$,
where $x$ and $C$ are respectively a letter and an admissible column. This
will be called `` the insertion of the letter $x$ in the admissible column $C$
'' and denoted by $x\rightarrow C$. Then we will be able to obtain $P(w)$ when
$w=\mathrm{w}(T)x$ with $x$ a letter and $T$ an orthogonal tableau. This will
be called `` the insertion of the letter $x$ in the orthogonal tableau $T$ ''
and denoted by $x\rightarrow T$. Our construction of $P$ will be recursive, in
the sense that if $P(u)=T$ and $x$ is a letter, then $P(ux)=x\rightarrow T$.

\subsubsection{Insertion of a letter in an admissible column\label{x_in_C}}

Consider a word $w=$\textrm{w(}$C)x$, where $x$ and $C$ are respectively a
letter and an admissible column of height $p$. When $w=\mathrm{w}(C^{x})$ is
the reading of a column $C^{x}$, we have:
\begin{align*}
x  &  \rightarrow C=C^{x}\text{ if }C^{x}\text{ is admissible or}\\
x  &  \rightarrow C=\widetilde{C^{x}}\text{ where }\widetilde{C^{x}}\text{ is
the column whose reading corresponds to }\widetilde{w}\text{ otherwise.}%
\end{align*}
Indeed, $x\rightarrow C$ must be an orthogonal tableau such that
$\mathrm{w}(x\rightarrow C)\equiv w$.

When $w$ is not a column word, by Lemma \ref{lem_phi_tens} the highest weight
vertex $w^{0}$ of $B(w)\,$may be written $w^{0}=v^{0}\,1$ where $v^{0}%
\in\{b_{\omega_{p}};p=1,...,n\}\cup\{b_{\overline{\omega}_{n}}\}.$ Then
$u^{0}=1\,v^{0}$ is the reading of an orthogonal tableau and $u^{0}\equiv
w^{0}$.\ So $u^{0}$ is the highest weight vertex of the connected component
containing $\mathrm{w}(x\rightarrow C)$. Moreover there exists a unique
sequence of highest weight vertices $w_{1}^{0},...,w_{p}^{0}$ such that
$w_{1}^{0}=w^{0},$ $w_{p}=u^{0}$ and for $i=2,...,p$ $w_{i}^{0}$ differs from
$w_{i-1}^{0}$ by applying one relation $R_{1}$ from left to right.\ This
implies that there exists a unique sequence of vertices $w_{1},...,w_{p}$ such
that $w_{1}=w$ and for $i=2,...,p-1$ $B(w_{i})=B(w_{i}^{0})$. Each $w_{i}$
differs from $w_{i-1}$ by applying one relation $R_{1},R_{2},R_{3}$ or $R_{4}$
from left to right.\ The word $w_{p}$ is the reading of an orthogonal tableau
and can be factorized as $w_{p}=v^{\prime}\,x^{\prime}$ where $v^{\prime
}=\mathrm{w}(C^{\prime})$ is a column word an $x^{\prime}$ a letter. We will
have $x\rightarrow C=C^{\prime}x^{\prime}$.

\begin{example}
\ \ \ \ 

\noindent Suppose $n=7$. Let $\mathrm{w}(C)=6700\bar{7}\bar{6}$ be an
admissible column word of type $B$.\ Choose $x=6.$ Then by applying relations
$R_{i}^{B}$ $i=1,...,4$ we obtain successively:

\noindent$6700\mathbf{\bar{7}\bar{6}6}\equiv670\mathbf{0\bar{7}7}\bar{7}%
\equiv67\mathbf{0\bar{7}7}0\bar{7}\equiv6\mathbf{7\bar{7}7}00\bar{7}%
\equiv\mathbf{6\bar{6}6}700\bar{7}\equiv\bar{5}56700\bar{7}$

\noindent Suppose $n=7$. Let $\mathrm{w}(C)=67\bar{7}7\bar{7}\bar{6}$ be an
admissible column word of type $D$.\ Choose $x=6.$ Then by applying relations
$R_{i}^{D}$ $i=1,...,4$ we obtain successively:

\noindent$67\bar{7}7\mathbf{\bar{7}\bar{6}6}\equiv67\bar{7}\mathbf{7\bar
{7}\bar{7}}7\equiv67\mathbf{\bar{7}\bar{6}6}\bar{7}7\equiv6\mathbf{7\bar
{7}\bar{7}}7\bar{7}7\equiv\mathbf{6\bar{6}6}\bar{7}7\bar{7}7\equiv\bar
{5}567\bar{7}7\bar{7}.$

Hence%

\[
6\rightarrow%
\begin{tabular}
[c]{|l|}\hline
$\mathtt{6}$\\\hline
$\mathtt{7}$\\\hline
$\mathtt{0}$\\\hline
$\mathtt{0}$\\\hline
$\mathtt{\bar{7}}$\\\hline
$\mathtt{\bar{6}}$\\\hline
\end{tabular}
=%
\begin{tabular}
[c]{|l|l}\hline
$\mathtt{5}$ & \multicolumn{1}{|l|}{$\mathtt{\bar{5}}$}\\\hline
$\mathtt{6}$ & \\\cline{1-1}%
$\mathtt{7}$ & \\\cline{1-1}%
$\mathtt{0}$ & \\\cline{1-1}%
$\mathtt{0}$ & \\\cline{1-1}%
$\mathtt{\bar{7}}$ & \\\cline{1-1}%
\end{tabular}
\text{ and }6\rightarrow%
\begin{tabular}
[c]{|l|}\hline
$\mathtt{6}$\\\hline
$\mathtt{7}$\\\hline
$\mathtt{\bar{7}}$\\\hline
$\mathtt{7}$\\\hline
$\mathtt{\bar{7}}$\\\hline
$\mathtt{\bar{6}}$\\\hline
\end{tabular}
=%
\begin{tabular}
[c]{|l|l}\hline
$\mathtt{5}$ & \multicolumn{1}{|l|}{$\mathtt{\bar{5}}$}\\\hline
$\mathtt{6}$ & \\\cline{1-1}%
$\mathtt{7}$ & \\\cline{1-1}%
$\mathtt{\bar{7}}$ & \\\cline{1-1}%
$\mathtt{7}$ & \\\cline{1-1}%
$\mathtt{\bar{7}}$ & \\\cline{1-1}%
\end{tabular}
.
\]
\end{example}

\subsubsection{Insertion of a letter in an orthogonal tableau}

Consider an orthogonal tableau $T=C_{1}C_{2}\cdot\cdot\cdot C_{r}$. We can
prove as in \cite{Lec} that the insertion $x\rightarrow T$ is characterized as follows:

\begin{itemize}
\item  If $\mathrm{w(}C_{1})\,x$ is an admissible column word , then
$x\rightarrow T=C_{1}^{x}C_{2}\cdot\cdot\cdot C_{r}$ where $C_{1}^{x}$ is the
column of reading $\mathrm{w(}C_{1})\,x.$

\item  If $\mathrm{w(}C_{1})\,x$ is a non admissible column word each strict
factor of which is admissible and such that $\widetilde{x\mathrm{w(}C_{1}%
)}=x_{1}\cdot\cdot\cdot x_{s},$ then $x\rightarrow T=x_{s}\rightarrow
(x_{s-1}\rightarrow(\cdot\cdot\cdot x_{1}\rightarrow T^{\prime}))$ where
$T^{\prime}=C_{2}\cdot\cdot\cdot C_{r}$. Moreover the insertion of
$x_{1},...,x_{s}$ in $T^{\prime}$ does not cause a new contraction.

\item  If $\mathrm{w(}C_{1})\,x$ is not a column word, the insertion of $x$ in
$C_{1}$ gives a column $C_{1}^{\prime}$ and a letter $x^{\prime}$ (with the
notation of \ref{x_in_C}). Then $x\rightarrow T=C_{1}^{\prime}(x^{\prime
}\rightarrow T^{\prime})$, that is, $x\rightarrow T$ is the tableau defined by
$C_{1}^{\prime}$ and the columns of $x^{\prime}\rightarrow T^{\prime}.$
\end{itemize}

Notice that the algorithm terminates because in the last two cases we are
reduced to the insertion of a letter in a tableau whose number of boxes is
strictly less than that of $T$. Finally for any vertex $w\in G_{n}$, we will
have:
\begin{align*}
P(w)  &  =
\begin{tabular}
[c]{|l|}\hline
$w$\\\hline
\end{tabular}
\text{ if }w\text{ is a letter,}\\
P(w)  &  =x\rightarrow P(u)\text{ if }w=ux\text{ with }u\text{ a word and
}x\text{ a letter.}%
\end{align*}

\subsection{Schensted-type Correspondences \label{sub_sec_Cor_in_Gn}}

In this section a bijection is established between words $w$ of length $l$ on
$\mathcal{B}_{n}\ $and pairs $(P^{B}(w),Q^{B}(w))$ where $P^{B}(w)$ is the
orthogonal tableau defined above and $Q^{B}(w)$ is an oscillating tableau of
type $B$. Similarly we obtain a bijection between words $w$ of length $l$ on
$\mathcal{D}_{n}\ $and pairs $(P^{D}(w),Q^{D}(w))$ where $P^{D}(w)$ is an
oscillating tableau of type $D$. For type $B$, such a one-to-one
correspondence has already been obtained by Sundaram \cite{SU} using another
definition of orthogonal tableaux and an appropriate insertion algorithm.
Unfortunately it is not known if this correspondence is compatible with a
monoid structure. Our bijection based on the previous insertion algorithm for
admissible orthogonal tableaux of type $B$ will be different from Sundaram's
one but compatible with the plactic relations defining $Pl(B_{n})$.

\begin{definition}
\label{def_tab_osci} \ \ \ \ 

\noindent An oscillating tableau $Q$ of type $B$ and length $l$ is a sequence
of Young diagrams $(Q_{1},...,Q_{l})$ whose columns have height $\leq n$ and
such that any two consecutive diagrams are equal or differ by exactly one box
(i.e. $Q_{k+1}=Q_{k}$, $Q_{k+1}/Q_{k}=(\,%
\begin{tabular}
[c]{|l|}\hline
\\\hline
\end{tabular}
\,)$ or $Q_{k}/Q_{k+1}=(\,%
\begin{tabular}
[c]{|l|}\hline
\\\hline
\end{tabular}
\,)$).

\noindent An oscillating tableau $Q$ of type $D$ and length $l$ is a sequence
$(Q_{1},...,Q_{l})$ of pairs $Q_{k}(O_{k},\varepsilon_{k})$ where $O_{k}$ is a
Young diagram whose columns have height $\leq n$ and $\varepsilon_{k}%
\in\{-,0,+\},$ satisfying for $k=1,...,l$

\begin{itemize}
\item $O_{k+1}/O_{k}=(\,%
\begin{tabular}
[c]{|l|}\hline
\\\hline
\end{tabular}
\,)$ or $O_{k}/O_{k+1}=(\,%
\begin{tabular}
[c]{|l|}\hline
\\\hline
\end{tabular}
\,)$,

\item $\varepsilon_{k+1}\neq0$ and $\varepsilon_{k}\neq0$ implies
$\varepsilon_{k+1}=\varepsilon_{k}$.

\item $\varepsilon_{k}=0$ if and only if $O_{k}$ has no columns of height $n.$
\end{itemize}
\end{definition}

Let $w=x_{1}\cdot\cdot\cdot x_{l}\ $be a word. The construction of $P(w)$
involves the construction of the $l$ orthogonal tableaux defined by
$P_{i}=P(x_{1}\cdot\cdot\cdot x_{i}).$ For $w\in\mathcal{B}_{n}^{\ast}$ (resp.
$w\in\mathcal{D}_{n}^{\ast}$) we denote by $Q^{B}(w)$ (resp. $Q^{D}(w) $) the
sequence of shapes of the orthogonal tableaux $P_{1},...,P_{l}$.

\begin{proposition}
\label{prop_Q(w)_oscill}$Q_{B}(w)$ and $Q_{D}(w)$ are respectively oscillating
tableaux of type $B$ and $D$.
\end{proposition}

\begin{proof}
Each $Q_{i}$ is the shape of an orthogonal tableau so it suffices to prove
that for any letter $x$ and any orthogonal tableau $T$, the shape of
$x\rightarrow T$ differs from the shape of $T$ by at most one box according to
Definition \ref{def_tab_osci}.

The highest weight vertex of the connected component containing $\mathrm{w(}%
T)x$ may be written $\mathrm{w(}T^{0})x^{0}$ where $T^{0}$ is an orthogonal
tableau. It follows from Lemma \ref{lem_coplactic} $\mathrm{(ii)}$ that
$\mathrm{w(}T)\longleftrightarrow\mathrm{w(}T^{0})$. So $\mathrm{wt}%
(\mathrm{w(}T^{0}))$ is given by the shape of $T$. Then the shape of
$x\rightarrow T$ is given by the coordinates of $\mathrm{wt}(\mathrm{w}%
(T^{0})x^{0})$ on the basis $(\omega_{1}^{B},...,\omega_{n}^{B})$ for type
$B$, on the base $(\omega_{1}^{D},...,\omega_{n}^{D})$ or $(\omega_{1}%
^{D},...,\omega_{n-1}^{D},\overline{\omega}_{n}^{D})$ for type $D$.

Suppose that $x\in\mathcal{B}_{n}^{\ast}$ and $T$ is orthogonal of type $B$.
Let $(\lambda_{1},...,\lambda_{n})$ be the coordinates of $\mathrm{wt}(T^{0})$
on the basis of the $\omega_{i}^{B}$'s. If $x^{0}=\overline{i}\succ0$ then
$\mathrm{wt}(x^{0})=\omega_{i-1}^{B}-\omega_{i}^{B}$. So $\lambda_{i}>0$ and
$\mathrm{wt}(\mathrm{w(}T^{0})x^{0})=(\lambda_{1},...,\lambda_{i-1}%
+1,\lambda_{i}-1,...,\lambda_{n-1}).$ Hence during the insertion of the letter
$x$ in $T$, a column of height $i$ (corresponding to the weight $\omega_{i}$)
is turned into a column of height $i-1$ (corresponding to the weight
$\omega_{i-1}$)$.$ So the shape of $x\rightarrow T$ is obtained by erasing one
box to the shape of $T$. If $x^{0}=i\prec0$, then we can prove by similar
arguments that the shape of $x\rightarrow T$ is obtained by adding one box to
the shape of $T$. When $x^{0}=0,$ $\mathrm{wt}(x^{0})=0$, so $\mathrm{wt}%
(\mathrm{w(}T^{0})x^{0})=\mathrm{wt}(\mathrm{w(}T^{0}))$. Hence the shapes of
$T$ and $x\rightarrow T$ are the same.

Suppose $x\in\mathcal{D}_{n}^{\ast}$ and $T$ orthogonal of type $D$. When
$\left|  x^{0}\right|  \neq n,$ the proof is the same as above. If $x^{0}=n,$
$\mathrm{wt}(x^{0})=\Lambda_{n}-\Lambda_{n-1}=\omega_{n}-\omega_{n-1}%
=\omega_{n-1}-\overline{\omega}_{n}.$ We have to consider three cases,
$\mathrm{(i)}$: $\varepsilon_{T}=-$; $\mathrm{(ii)}$: $\varepsilon_{T}=0$ and
$\mathrm{(iii)}$: $\varepsilon_{T}=+$. Denote by $(\lambda_{1},...,\lambda
_{n})$ the positive decomposition of $\mathrm{wt}(\mathrm{w(}T^{0}))$ on the
basis $(\omega_{1}^{D},...,\omega_{n}^{D})$ or on the basis $(\omega_{1}%
^{D},...,\overline{\omega}_{n}^{D})$.

\noindent In the first case, $\lambda_{n}>0$ and the positive decomposition of
$\mathrm{wt}(x^{0}\mathrm{w(}T^{0}))$ on the basis $(\omega_{1}^{D}%
,...,\overline{\omega}_{n}^{D})$ is $(\lambda_{1},...,\lambda_{n-2}%
,\lambda_{n-1}+1,\lambda_{n}-1)$. It means that during the insertion of $x$ in
$T$ a column of height $n$ (corresponding to $\overline{\omega}_{n}$) is
turned into a column of height $n-1$ (corresponding to $\omega_{n-1}$).
Moreover $\varepsilon_{x\rightarrow T}=\varepsilon_{T}$ if $\lambda_{n}>1$ and
$\varepsilon_{x\rightarrow T}=0$ otherwise.

\noindent In the second case, $\lambda_{n-1}>0,$ $\lambda_{n}=0$ and the
positive decomposition of $\mathrm{wt}(x^{0}\mathrm{w(}T^{0}))$ on the base
$(\omega_{1}^{D},...,\omega_{n}^{D})$ is $(\lambda_{1},\lambda_{2}%
,...,\lambda_{n-1}-1,1).$ It means that during the insertion of $x$ in $T$ a
column of height $n-1$ (corresponding to $\omega_{n-1}$) is turned into a
column of height $n$ (corresponding to $\omega_{n}$). Moreover $\varepsilon
_{x\rightarrow T}=+$.

\noindent In the last case, $\lambda_{n-1}>0,$ $\lambda_{n}>0$ and the
positive decomposition of $\mathrm{wt}(x^{0}\mathrm{w}(T^{0}))$ on
$(\omega_{1}^{D},...,\omega_{n}^{D})$ is $(\lambda_{1},\lambda_{2}%
,...,\lambda_{n-1}-1,\lambda_{n}+1)$. It means that during the insertion of
$x$ in $T$ a column of height $n-1$ (corresponding to $\omega_{n-1}$) is
turned into a column of height $n$ (corresponding to $\omega_{n}$). Moreover
$\varepsilon_{x\rightarrow T}=\varepsilon_{T}$.

When $x^{0}=\overline{n},$ the proof is similar.
\end{proof}

\begin{theorem}
\label{th_Q_symbol}For any vertices $w_{1}$ and $w_{2}$ of $G_{n}$:
\[
w_{1}\longleftrightarrow w_{2}\Leftrightarrow Q(w_{1})=Q(w_{2})\text{.}%
\]
\end{theorem}

\begin{proof}
The proof is analogous to that of Proposition 5.2.1 in \cite{Lec}.
\end{proof}

\begin{corollary}
Let $\mathcal{B}_{n,l}^{\ast}$ and $\mathcal{O}_{l}^{B}$ (resp. $\mathcal{D}%
_{n,l}^{\ast}$ and $\mathcal{O}_{l}^{D}$) be the set of words of length $l$ on
$\mathcal{B}_{n}$ (resp. $\mathcal{D}_{n}$) and the set of pairs $(P,Q)$ where
$P$ is an orthogonal tableau of type $B$ (resp.\ $D$) and $Q$ an oscillating
tableau of type $B$ (resp. $D$) and length $l$ such that $P$ has shape $Q_{l}$
($Q_{l}$ is the last shape of $Q$). Then the maps:
\[
\mathcal{%
\begin{tabular}
[c]{l}%
$\Psi^{B}:B_{n,l}^{\ast}\rightarrow O_{l}^{B}$\\
$w\mapsto(P^{B}(w),Q^{B}(w))$%
\end{tabular}
}\text{ and }\mathcal{%
\begin{tabular}
[c]{l}%
$\Psi^{D}:\mathcal{D}_{n,l}^{\ast}\rightarrow\mathcal{O}_{l}^{D}$\\
$w\mapsto(P^{D}(w),Q^{D}(w))$%
\end{tabular}
}%
\]
are bijections.
\end{corollary}

\begin{proof}
For type $\Psi^{B}$ the proof is analogous to that of Theorem 5.2.2 in
\cite{Lec}. By Theorems \ref{th_Psymbol} and \ref{th_Q_symbol}, we obtain that
$\Psi^{D}$ is injective.\ Consider an oscillating tableau $Q$ of length $l$
and type $D$.\ Set $x_{1}=1$ and for $i=2,...,l$

\noindent- $x_{i}=k$ if $O_{i}$ differs from $O_{i-1}$ by adding a box in row
$k$ of height $<n,$

\noindent- $x_{i}=\overline{k}$ if $Q_{i}$ differs from $Q_{i-1}$ by removing
a box in row $k.$of\ height $<n,$

\noindent- $x_{i}=n$ if $O_{i}$ differs from $O_{i-1}$ by adding a box in row
$n$ and $\varepsilon_{i}=+,$

\noindent- $x_{i}=\overline{n}$ if $Q_{i}$ differs from $Q_{i-1}$ by adding a
box in row $n$ and $\varepsilon_{i}=-,$

\noindent- $x_{i}=\overline{n}$ if $O_{i}$ differs from $O_{i-1}$ by removing
a box in row $n$ and $\varepsilon_{i}=+,$

\noindent- $x_{i}=n$ if $O_{i}$ differs from $O_{i-1}$ by removing a box in
row $n$ and $\varepsilon_{i}=-,$

\noindent- Consider $w_{Q}=x_{l}\cdot\cdot\cdot x_{2}1$.\ Then $Q(w_{Q}%
)=Q$.\ By Theorem \ref{TH_KN}, the image of $B(w_{Q})$ by $\Psi^{D}$ consists
in the pairs $(P,Q)$ where $P$ is a symplectic tableau of shape $Q_{l}$.\ We
deduce immediately that $\Psi$ is surjective.
\end{proof}

\subsection{\label{subsecJDT}Jeu de Taquin for type B}

In \cite{SH}, J T Sheats has developed a sliding algorithm for type $C$ acting
on the skew admissible symplectic tableaux. This algorithm is analogous to the
classical Jeu de Taquin of Lascoux and Sch\"{u}tzenberger for type A
\cite{LS}. Each inner corner of the skew tableau considered is turned into an
outside corner by applying vertical and horizontal moves. We have shown in
\cite{Lec} how to extend it to take into account the contraction relation of
the plactic monoid $Pl(C_{n})$ (analogous to $Pl(B_{n})$ and $Pl(D_{n})$ for
type $C$). Then we have proved that the tableau obtained does not depend on
the way the inner corners disappear. In this section we propose a sliding
algorithm for type $B$. The main idea is that the split form of any skew
orthogonal tableau $T$ of type $B$ may be regarded as a symplectic skew tableau.

Set $\mathcal{C}_{n}=\{1\prec\cdot\cdot\cdot\prec n\prec\overline{n}\prec
\cdot\cdot\cdot\prec\overline{1}\}\subset\mathcal{B}_{n}$. The symplectic
tableaux are, for type $C,$ the combinatorial objects analogous to the
orthogonal tableaux. They can be regarded as orthogonal tableaux of type $B$
on the alphabet $\mathcal{C}_{n}$ instead of $\mathcal{B}_{n}$. The plactic
monoid $Pl(C_{n})$ is the quotient of the free monoid $\mathcal{C}_{n}^{\ast}
$ by relations $R_{1}^{B},$ $R_{2}^{B}$ and $R_{5}^{B}$. We denote by
$\overset{C}{\equiv}$ the congruence relation in $Pl(C_{n})$. Then for $w_{1}$
and $w_{2}$ two words of $\mathcal{C}_{n}^{\ast}$ we have:
\[
w_{1}\overset{C}{\equiv}w_{2}\Longrightarrow w_{1}\overset{B}{\equiv}%
w_{2}\text{.}%
\]
A skew orthogonal tableau of type $B$ is a skew Young diagram filled by
letters of $\mathcal{B}_{n}$ whose columns are admissible of type $B$ and such
that the rows of its split form (obtained by splitting its columns) are weakly
increasing from left to right. Skew orthogonal tableaux are the combinatorial
objects analogous to the admissible skew tableaux introduced by Sheats in
\cite{SH} for type $C$. Note that two different skew tableaux may have the
same reading.

\begin{example}
For $n=3,$

\noindent$T=%
\begin{tabular}
[c]{lll}\cline{3-3}%
&  & \multicolumn{1}{|l|}{$\mathtt{2}$}\\\cline{2-3}%
& \multicolumn{1}{|l}{$\mathtt{3}$} & \multicolumn{1}{|l|}{$\mathtt{0}$%
}\\\hline
\multicolumn{1}{|l}{$\mathtt{0}$} & \multicolumn{1}{|l}{$\mathtt{\bar{3}}$} &
\multicolumn{1}{|l|}{$\mathtt{\bar{1}}$}\\\hline
\multicolumn{1}{|l}{$\mathtt{0}$} & \multicolumn{1}{|l}{} & \\\cline{1-1}%
\end{tabular}
$ is a skew orthogonal tableau of type $B$ because $\mathrm{spl}%
(T)=\vspace{0.2cm}%
\begin{tabular}
[c]{llllll}\cline{5-6}%
&  &  &  & \multicolumn{1}{|l}{$\mathtt{2}$} &
\multicolumn{1}{|l|}{$\mathtt{2}$}\\\cline{3-6}%
&  & \multicolumn{1}{|l}{$\mathtt{2}$} & \multicolumn{1}{|l}{$\mathtt{3}$} &
\multicolumn{1}{|l}{$\mathtt{3}$} & \multicolumn{1}{|l|}{$\mathtt{\bar{3}}$%
}\\\hline
\multicolumn{1}{|l}{$\mathtt{2}$} & \multicolumn{1}{|l}{$\mathtt{\bar{3}}$} &
\multicolumn{1}{|l}{$\mathtt{\bar{3}}$} & \multicolumn{1}{|l}{$\mathtt{\bar
{2}}$} & \multicolumn{1}{|l}{$\mathtt{\bar{1}}$} &
\multicolumn{1}{|l|}{$\mathtt{\bar{1}}$}\\\hline
\multicolumn{1}{|l}{$\mathtt{3}$} & \multicolumn{1}{|l}{$\mathtt{\bar{2}}$} &
\multicolumn{1}{|l}{} &  &  & \\\cline{1-2}%
\end{tabular}
$.
\end{example}

The relation $0\overline{n}n\equiv\overline{n}n0$ has no natural
interpretation in terms of horizontal or vertical slidings in skew orthogonal
tableaux. To overcome this problem we are going to work on the split form of
the skew tableaux instead of the skew tableaux themselves that is , we are
going to obtain a Jeu de Taquin for type $B$ by applying the symplectic Jeu de
Taquin on the split form of the skew orthogonal tableaux of type $B$.

\begin{lemma}
Let $T$ and $T^{\prime}$ be two skew orthogonal tableaux of type $B$. Then:
\[
\mathrm{w}(T)\overset{B}{\equiv}\mathrm{w}(T^{\prime})\Longleftrightarrow
\mathrm{w}[\mathrm{spl}(T)]\overset{B}{\equiv}\mathrm{w}[\mathrm{spl}%
(T^{\prime})].
\]
\end{lemma}

\begin{proof}
We can write $\mathrm{w}(T)=\mathrm{w}(C_{1})\cdot\cdot\cdot\mathrm{w}(C_{r})$
and $\mathrm{w}(T^{\prime})=\mathrm{w}(C_{1}^{\prime})\cdot\cdot
\cdot\mathrm{w}(C_{s}^{\prime})$ where $C_{k}$ and $C_{k}^{\prime},$
$k=1,...,r$ are admissible columns. All the vertices $w\in B(\mathrm{w}(T))$
and $w^{\prime}\in B(\mathrm{w}(T^{\prime}))$ can be respectively written on
the form $w=c_{r}\cdot\cdot\cdot c_{1}$ and $w^{\prime}=c_{s}^{\prime}%
\cdot\cdot\cdot c_{1}^{\prime}$ where $c_{i},$ $i=1,..,r$ and $c_{j}^{\prime}%
$, $j=1,...,s$ are readings of admissible columns of type $B$.\ Consider the
maps:
\[
\theta_{2}:\left\{
\begin{tabular}
[c]{c}%
$B(\mathrm{w}(T))\rightarrow B(\mathrm{spl}(\mathrm{w}(T))$\\
$w=c_{r}\cdot\cdot\cdot c_{1}\longmapsto S_{2}(c_{r})\cdot\cdot\cdot
S_{2}(c_{1})$%
\end{tabular}
\right.  \text{ and }\theta_{2}^{\prime}:\left\{
\begin{tabular}
[c]{c}%
$B(\mathrm{w}(T^{\prime}))\rightarrow B(\mathrm{spl}(\mathrm{w}(T))$\\
$w^{\prime}=c_{s}^{\prime}\cdot\cdot\cdot c_{1}\longmapsto S_{2}(c_{r}%
^{\prime})\cdot\cdot\cdot S_{2}(c_{1}^{\prime})$%
\end{tabular}
\right.
\]
where $S_{2}$ is the map defined in Proposition \ref{prop_imag_S2}.\ We have
$\mathrm{w}[\mathrm{spl}(T)]=\theta_{2}(\mathrm{w}(T))$ and $\mathrm{w}%
[\mathrm{spl}(T^{\prime})]=\theta_{2}^{\prime}(\mathrm{w}(T^{\prime})).$ By
using \ Corollary \ref{cor_strech} we obtain
\[
\mathrm{w}(T)\overset{B}{\equiv}\mathrm{w}(T^{\prime})\Longleftrightarrow
\mathrm{w}(T)\overset{B}{\sim}\mathrm{w}(T^{\prime})\Longleftrightarrow
\mathrm{w}[\mathrm{spl}(T)]\overset{B}{\sim}\mathrm{w}[\mathrm{spl}(T^{\prime
})]\Longleftrightarrow\mathrm{w}[\mathrm{spl}(T)]\overset{B}{\equiv}%
\mathrm{w}[\mathrm{spl}(T^{\prime})].
\]
\end{proof}

If $T$ is a skew orthogonal tableau of type $B$ with $r$ columns, then
$\mathrm{spl}(T)$ is a symplectic skew tableau with $2r$ columns. We can apply
the symplectic Jeu de taquin to $\mathrm{spl}(T)$ to obtain a symplectic
tableau $\mathrm{spl}(T)^{\prime}$. We will have $\mathrm{w}[\mathrm{spl}%
(T)^{\prime}]\overset{C}{\equiv}\mathrm{w}[\mathrm{spl}(T)]$ so $\mathrm{w}%
[\mathrm{spl}(T)^{\prime}]\overset{B}{\equiv}\mathrm{w}[\mathrm{spl}(T)]$.

\begin{proposition}
$\mathrm{spl}(T)^{\prime}$ is the split form of the orthogonal tableau
$P^{B}(T)$.
\end{proposition}

\begin{proof}
It follows from $\mathrm{w}(T)\overset{B}{\equiv}\mathrm{w}(P_{B}(T))$ and the
lemma above that $\mathrm{w}[\mathrm{spl}(T)]\overset{B}{\equiv}%
\mathrm{w}[\mathrm{spl}(P^{B}(T))].$ So we obtain $\mathrm{w}[\mathrm{spl}%
(T)^{\prime}]\overset{B}{\equiv}\mathrm{w}[\mathrm{spl}(P^{B}(T))]$. But
$\mathrm{spl}(T^{\prime})$ and $\mathrm{spl}(P^{B}(T))$ are orthogonal
tableaux, hence $\mathrm{spl}(T)^{\prime}=\mathrm{spl}(P^{B}(T)).$
\end{proof}

\noindent The columns of the split form of a skew orthogonal tableau $T$ of
type $B$ contain no letters $0$ and no pairs of letters $(x,\overline{x})$
with $x\preceq n.$ In this particular case most of the elementary steps of the
symplectic Jeu de Taquin applied on $T$ are simple slidings identical to those
of the original Jeu de Taquin of Lascoux and Sch\"{u}tzenberger (that is
complications of the symplectic Jeu de taquin are not needed in these slidings).

\begin{example}
From $\mathrm{spl}\left(
\begin{tabular}
[c]{l|l|l|}\cline{2-3}\cline{3-3}%
& $\mathtt{1}$ & $\mathtt{2}$\\\hline
\multicolumn{1}{|l|}{$\mathtt{1}$} & $\mathtt{0}$ & $\mathtt{\bar{3}}$\\\hline
\multicolumn{1}{|l|}{$\mathtt{3}$} & $\mathtt{\bar{3}}$ & $\mathtt{\bar{2}}%
$\\\hline
\end{tabular}
\right)  =%
\begin{tabular}
[c]{|l|l|l|l|l|l|}\hline
$\mathtt{\ast}$ & $\mathtt{\ast}$ & $\mathtt{1}$ & $\mathtt{1}$ & $\mathtt{1}$%
& $\mathtt{2}$\\\hline
$\mathtt{1}$ & $\mathtt{1}$ & $\mathtt{2}$ & $\mathtt{\bar{3}}$ &
$\mathtt{\bar{3}}$ & $\mathtt{\bar{3}}$\\\hline
$\mathtt{3}$ & $\mathtt{3}$ & $\mathtt{\bar{3}}$ & $\mathtt{\bar{2}}$ &
$\mathtt{\bar{2}}$ & $\mathtt{\bar{1}}$\\\hline
\end{tabular}
$ we compute successively:\vspace{0.5cm}

$%
\begin{tabular}
[c]{|l|l|l|l|l|l|}\hline
$\mathtt{\ast}$ & $\mathtt{1}$ & $\mathtt{1}$ & $\mathtt{1}$ & $\mathtt{1}$ &
$\mathtt{2}$\\\hline
$\mathtt{1}$ & $\mathtt{2}$ & $\mathtt{\ast}$ & $\mathtt{\bar{3}}$ &
$\mathtt{\bar{3}}$ & $\mathtt{\bar{3}}$\\\hline
$\mathtt{3}$ & $\mathtt{3}$ & $\mathtt{\bar{3}}$ & $\mathtt{\bar{2}}$ &
$\mathtt{\bar{2}}$ & $\mathtt{\bar{1}}$\\\hline
\end{tabular}
,$ \vspace{0.5cm}$%
\begin{tabular}
[c]{|l|l|l|l|l|l|}\hline
$\mathtt{\ast}$ & $\mathtt{1}$ & $\mathtt{1}$ & $\mathtt{1}$ & $\mathtt{1}$ &
$\mathtt{2}$\\\hline
$\mathtt{1}$ & $\mathtt{2}$ & $\mathtt{\bar{3}}$ & $\mathtt{\bar{3}}$ &
$\mathtt{\bar{3}}$ & $\mathtt{\bar{3}}$\\\hline
$\mathtt{3}$ & $\mathtt{3}$ & $\mathtt{\bar{2}}$ & $\mathtt{\ast}$ &
$\mathtt{\bar{2}}$ & $\mathtt{\bar{1}}$\\\hline
\end{tabular}
,$ $%
\begin{tabular}
[c]{|l|l|l|l|l|l|}\hline
$\mathtt{\ast}$ & $\mathtt{1}$ & $\mathtt{1}$ & $\mathtt{1}$ & $\mathtt{1}$ &
$\mathtt{2}$\\\hline
$\mathtt{1}$ & $\mathtt{2}$ & $\mathtt{\bar{3}}$ & $\mathtt{\bar{3}}$ &
$\mathtt{\bar{3}}$ & $\mathtt{\bar{3}}$\\\hline
$\mathtt{3}$ & $\mathtt{3}$ & $\mathtt{\bar{2}}$ & $\mathtt{\bar{2}}$ &
$\mathtt{\ast}$ & $\mathtt{\bar{1}}$\\\hline
\end{tabular}
$

$%
\begin{tabular}
[c]{|l|l|l|l|l|l|}\hline
$\mathtt{\ast}$ & $\mathtt{1}$ & $\mathtt{1}$ & $\mathtt{1}$ & $\mathtt{2}$ &
$\mathtt{2}$\\\hline
$\mathtt{1}$ & $\mathtt{2}$ & $\mathtt{\bar{3}}$ & $\mathtt{\bar{3}}$ &
$\mathtt{\bar{3}}$ & $\mathtt{\bar{3}}$\\\hline
$\mathtt{3}$ & $\mathtt{3}$ & $\mathtt{\bar{2}}$ & $\mathtt{\bar{2}}$ &
$\mathtt{\bar{2}}$ & $\mathtt{\ast}$\\\hline
\end{tabular}
,$ $%
\begin{tabular}
[c]{|l|l|l|l|l|l|}\hline
$\mathtt{1}$ & $\mathtt{1}$ & $\mathtt{1}$ & $\mathtt{1}$ & $\mathtt{2}$ &
$\mathtt{2}$\\\hline
$\mathtt{2}$ & $\mathtt{\ast}$ & $\mathtt{\bar{3}}$ & $\mathtt{\bar{3}}$ &
$\mathtt{\bar{3}}$ & $\mathtt{\bar{3}}$\\\hline
$\mathtt{3}$ & $\mathtt{3}$ & $\mathtt{\bar{2}}$ & $\mathtt{\bar{2}}$ &
$\mathtt{\bar{2}}$ & $\mathtt{\ast}$\\\hline
\end{tabular}
,$ $%
\begin{tabular}
[c]{|l|l|l|l|l|l|}\hline
$\mathtt{1}$ & $\mathtt{1}$ & $\mathtt{1}$ & $\mathtt{1}$ & $\mathtt{2}$ &
$\mathtt{2}$\\\hline
$\mathtt{2}$ & $\mathtt{3}$ & $\mathtt{\bar{3}}$ & $\mathtt{\bar{3}}$ &
$\mathtt{\bar{3}}$ & $\mathtt{\bar{3}}$\\\hline
$\mathtt{3}$ & $\mathtt{\bar{2}}$ & $\mathtt{\ast}$ & $\mathtt{\bar{2}}$ &
$\mathtt{\bar{2}}$ & $\mathtt{\ast}$\\\hline
\end{tabular}
,\vspace{0.5cm}$

$%
\begin{tabular}
[c]{|l|l|l|l|l|l|}\hline
$\mathtt{1}$ & $\mathtt{1}$ & $\mathtt{1}$ & $\mathtt{1}$ & $\mathtt{2}$ &
$\mathtt{2}$\\\hline
$\mathtt{2}$ & $\mathtt{3}$ & $\mathtt{\bar{3}}$ & $\mathtt{\bar{3}}$ &
$\mathtt{\bar{3}}$ & $\mathtt{\bar{3}}$\\\hline
$\mathtt{3}$ & $\mathtt{\bar{2}}$ & $\mathtt{\bar{2}}$ & $\mathtt{\ast}$ &
$\mathtt{\bar{2}}$ & $\mathtt{\ast}$\\\hline
\end{tabular}
,$ $%
\begin{tabular}
[c]{|l|l|l|l|l|l|}\hline
$\mathtt{1}$ & $\mathtt{1}$ & $\mathtt{1}$ & $\mathtt{1}$ & $\mathtt{2}$ &
$\mathtt{2}$\\\hline
$\mathtt{2}$ & $\mathtt{3}$ & $\mathtt{\bar{3}}$ & $\mathtt{\bar{3}}$ &
$\mathtt{\bar{3}}$ & $\mathtt{\bar{3}}$\\\hline
$\mathtt{3}$ & $\mathtt{\bar{2}}$ & $\mathtt{\bar{2}}$ & $\mathtt{\bar{2}}$ &
$\mathtt{\ast}$ & $\mathtt{\ast}$\\\hline
\end{tabular}
=\mathrm{spl}\left(
\begin{tabular}
[c]{|l|l|l}\hline
$\mathtt{1}$ & $\mathtt{1}$ & \multicolumn{1}{|l|}{$\mathtt{2}$}\\\hline
$\mathtt{3}$ & $\mathtt{\bar{3}}$ & \multicolumn{1}{|l|}{$\mathtt{\bar{3}}$%
}\\\hline
$\mathtt{0}$ & $\mathtt{\bar{2}}$ & \\\cline{1-2}%
\end{tabular}
\right)  .$ Note that the sliding applied in the fourth duplicated tableau
above is the unique sliding which is not identical to an original Jeu de
taquin step.
\end{example}

The split form of a skew orthogonal tableau of type $D$ (defined in the same
way than for type $B$) is still a symplectic skew tableau. But
\[
w_{1}\overset{C}{\equiv}w_{2}\not \Longrightarrow w_{1}\overset{D}{\equiv
}w_{2}%
\]
so we can not use the same idea to obtain an Jeu de Taquin for type $D$.
Moreover the examples (computed by using $P^{D}$ with $n=3$)
\[%
\begin{tabular}
[c]{|l|l|}\hline
1 & $\mathtt{3}$\\\hline
$\mathtt{\bar{3}}$ & $\mathtt{\bar{2}}$\\\hline
$\mathtt{\ast}$ & $\mathtt{\bar{1}}$\\\hline
\end{tabular}
\equiv%
\begin{tabular}
[c]{|l|l}\hline
$\mathtt{2}$ & \multicolumn{1}{|l|}{$\mathtt{3}$}\\\hline
$\mathtt{\bar{3}}$ & \multicolumn{1}{|l|}{$\mathtt{\bar{2}}$}\\\hline
$\mathtt{\bar{2}}$ & \\\cline{1-1}%
\end{tabular}
\text{ and }
\begin{tabular}
[c]{|l|l|}\hline
1 & $\mathtt{\bar{3}}$\\\hline
$\mathtt{\bar{3}}$ & $\mathtt{\bar{2}}$\\\hline
$\mathtt{\ast}$ & $\mathtt{\bar{1}}$\\\hline
\end{tabular}
\equiv%
\begin{tabular}
[c]{|l|l}\hline
$\mathtt{\bar{3}}$ & \multicolumn{1}{|l|}{$\mathtt{\bar{3}}$}\\\hline
$\mathtt{3}$ & \multicolumn{1}{|l|}{$\mathtt{\bar{2}}$}\\\hline
$\mathtt{\bar{3}}$ & \\\cline{1-1}%
\end{tabular}
\]
show that it is not enough to know what letter $x$ slides from the second
column $C_{2}$ to the first $C_{1}$ to be able to compute an horizontal
sliding. Indeed the result depends on the whole column $C_{2}$. Thus, to give
a combinatorial description of a sliding algorithm for type $D$ would probably
be very complicated.

\section{Plactic monoid for $\frak{G}_{n}\label{last_part}$}

Write $\frak{G}_{n}^{B}$ and $\frak{G}_{n}^{D}$ for the crystal graphs of the
direct sums
\[
\underset{l\geq0}{\bigoplus}(V(\Lambda_{1}^{B})%
{\textstyle\bigoplus}
V(\Lambda_{n}^{B}))^{\otimes l}\text{ and}\underset{l\geq0}{\text{ }\bigoplus
}(V(\Lambda_{1}^{D})%
{\textstyle\bigoplus}
V(\Lambda_{n}^{D})%
{\textstyle\bigoplus}
V(\Lambda_{n-1}^{D}))^{\otimes l}.
\]
We call $\frak{B}_{n}=\mathcal{B}_{n}\cup SP_{n}$ and $\frak{D}_{n}%
=\mathcal{D}_{n}\cup SP_{n}$ the sets of generalized letters of type $B$ and
$D$. Then we identify the vertices of $\frak{G}_{n}^{B}$ and $\frak{G}_{n}%
^{D}$ respectively with the words of the free monoid $\frak{B}_{n}^{\ast}$ and
$\frak{D}_{n}^{\ast}$. If $w$ is a vertex of $\frak{G}_{n},$ we write
$\mathrm{wt}(w)$ for the weight of $w$. The spin representations are
minuscule, hence every spin column is determined by its weight.

\noindent We can extend the Definition \ref{def_sam_plac} to vertices of
$\frak{G}_{n}$. Consider two vertices $\frak{b}_{1}$ and $\frak{b}_{2}$ of
$\frak{G}_{n}^{B}$ (resp. $\frak{G}_{n}^{D}$). We write $\frak{b}_{1}%
\overset{B}{\sim}\frak{b}_{2}$ (resp. $\frak{b}_{1}\overset{D}{\sim}%
\frak{b}_{2}$) when these vertices occur at the same place in two isomorphic
connected components of $\frak{G}_{n}^{B}$ (resp. $\frak{G}_{n}^{D}$). Our aim
is now to extend the results of Section \ref{subsec_monoids} to the vertices
of $\frak{G}_{n}$.

\subsection{Tensor products of spin representations}

Write $B(0)$ for the connected component of $\frak{G}_{n}$ containing only the
empty word.\ Let $\frak{C}_{0}$ be the spin column containing only barred
letters.\ For $p=1,...,n$, denote by $\frak{C}_{p}$ the spin column containing
exactly the unbarred letters $x\preceq p$. For any admissible column $C$, set
$\left|  C\right|  =\{x\preceq n,$ $x\in lC$ or $\overline{x}\in
lC\}=\{x\preceq n,$ $x\in rC$ or $\overline{x}\in rC\}.$

\begin{lemma}
\label{Lem_S_B} \ \ \ 

\begin{enumerate}
\item  There exists a unique crystal isomorphism $S^{B}$%
\[
B(0)\cup B(v_{\omega_{n}^{B}})\cup\left(  \overset{n-1}{\underset{i=1}{\cup}%
}B(v_{\omega_{i}^{B}})\right)  \overset{S^{B}}{\rightarrow}B(v_{\Lambda
_{n}^{B}})^{\otimes2}.
\]

\item  Let $w$ be the reading of an admissible column $C$ of type $B$.\ Write

- $l\frak{C}$ for the spin column of height $n$ obtained by adding to $lC$ the
barred letters $\overline{x}$ such that $x\notin\left|  C\right|  $,

- $r\frak{C}$ for the spin column of height $n$ obtained by adding to $rC$ the
unbarred letters $x$ such that $x\notin\left|  C\right|  $.
\end{enumerate}

\noindent Then
\[
S^{B}(w)=r\frak{C}\otimes l\frak{C.}%
\]
\end{lemma}

\begin{proof}
$1:$ From Lemma \ref{lem_phi_tens} we obtain that the highest weight vertices
of $B(v_{\Lambda_{n}^{B}})^{\otimes2}$ are the vertices $v_{p}^{B}%
=\frak{C}_{n}\otimes\frak{C}_{p}$ with $p=0,...,n$. We have $\mathrm{wt}%
(v_{p}^{B})=\omega_{p}^{B}$ for $p=1,...,n$ and $\mathrm{wt}(v_{0}^{D}%
)=0.\;$Hence $S^{B}$ is the crystal isomorphism which sends $B(v_{\omega
_{p}^{B}})$ on $B(v_{p}^{B})$ for $p=1,...,n$ and $B(0)$ on $B(v_{0}^{B}).$

$2:$When $w=v_{\omega_{p}^{B}}$, the equality $S^{B}(w)=r\frak{C}\otimes
l\frak{C}$ is true. Consider $w\in B(v_{\omega_{p}^{B}})$ and $i=1,...,n$ such
that $w^{\prime}=\widetilde{f}_{i}(w)\neq0$.\ Write $w=\mathrm{w}(C)$ and
$w^{\prime}=\mathrm{w}(C^{\prime})$ where $C$ and $C^{\prime}$ are two
admissible columns of height $p$.\ The lemma will be proved if we show the
implication
\[
S^{B}(w)=r\frak{C}\otimes l\frak{C\Longrightarrow} S^{B}(w^{\prime}%
)=r\frak{C}^{\prime}\otimes l\frak{C}^{\prime}%
\]
where $r\frak{C}^{\prime}$ and $l\frak{C}^{\prime}$ are defined from
$C^{\prime}$ in the same manner than $r\frak{C}$ and $l\frak{C}$ from
$C$.\ This is equivalent to
\begin{equation}
\widetilde{f}_{i}(r\frak{C}\otimes l\frak{C)=}r\frak{C}^{\prime}\otimes
l\frak{C}^{\prime}. \label{f_S}%
\end{equation}
Suppose $i\neq n$. Set $E_{i}=\{i,i+1,\overline{i+1},\overline{i}\}$.

$\mathrm{(i):}$ If $\{i,i+1\}\subset\left|  C\right|  $, $lC$ and $l\frak{C}$
coincide on $E_{i}$.\ Similarly $rC$ and $r\frak{C}$, $lC^{\prime}$ and
$l\frak{C}^{\prime}$, $rC^{\prime}$ and $l\frak{C}$' coincide on $E_{i}$. By
Proposition \ref{prop_imag_S2}, we know that
\[
\widetilde{f}_{i}^{2}(rC\otimes lC\frak{)=}rC^{\prime}\otimes lC^{\prime}.
\]
The action of $\widetilde{f}_{i}^{2}$ on $rC\otimes lC$ is analogous to that
of $\widetilde{f}_{i}$ on $r\frak{C}\otimes l\frak{C}$. It means that
$\widetilde{f}_{i}$ changes a pair $(i,\overline{i+1})$ of $r\frak{C}$ (resp
$l\frak{C})$ into a pair $(i+1,\overline{i})$ if and only if $\widetilde
{f}_{i}^{2}$ changes a pair $(i,\overline{i+1})$ of $rC$ (resp.\ $lC)$ into a
pair $(i+1,\overline{i})$.\ So (\ref{f_S}) is true because only the letters of
$E_{i}$ may be modified when we apply $\widetilde{f}_{i}.$

$\mathrm{(ii):}$ If $\{i,i+1\}\cap\left|  C\right|  =\{i+1\},$ we have
$[rC]_{i}=[lC]_{i}=\overline{i+1}$ with the notation of the proof of
Proposition \ref{prop_imag_S2}.\ Then $r\frak{C}\cap E_{i}=\{\overline
{i+1},i\}$ and $l\frak{C}\cap E_{i}=\{\overline{i+1},\overline{i}%
\}$.\ Moreover $[C^{\prime}]_{i}=\overline{i}$, $r\frak{C}^{\prime}\cap
E_{i}=\{\overline{i},i+1\}$ and $l\frak{C}^{\prime}\cap E_{i}=\{\overline
{i+1},\overline{i}\}$.\ Hence $\widetilde{f}_{i}(r\frak{C}\otimes l\frak{C)}$
and $r\frak{C}^{\prime}\otimes l\frak{C}^{\prime}$ coincide on $E_{i}$.\ So
they are equal because $\widetilde{f}_{i}$ does not modify the letters
$x\notin E_{i}$.

$\mathrm{(iii):}$ If $\{i,i+1\}\cap\left|  C\right|  =\{i\},$ the proof is
analogous to case $\mathrm{(ii).}$

Suppose $i=n$. Set $E_{n}=\{n,\overline{n}\}$.\ Then $n\in\left|  C\right|  $
because $\widetilde{f}_{i}(w)\neq0$.\ We obtain (\ref{f_S}) by using similar
arguments to those of $\mathrm{(i)}$.
\end{proof}

\begin{lemma}
\label{Lem_S_D} \ \ \ 

\begin{enumerate}
\item  There exists two crystal isomorphisms $S_{n}^{D}$ and $S_{n-1}^{D}$%
\begin{align*}
&  B(0)\cup B(v_{\omega_{n}^{D}})\cup\left(  \overset{n-1}{\underset{i=1}%
{\cup}}B(v_{\omega_{i}^{D}})\right)  \overset{S_{n}^{D}}{\rightarrow
}B(v_{\Lambda_{n}^{D}})\otimes(B(v_{\Lambda_{n}^{D}})\cup B(v_{\Lambda
_{n-1}^{D}})),\\
&  B(0)\cup B(v_{\overline{\omega}_{n}^{D}})\cup\left(  \overset
{n-1}{\underset{i=1}{\cup}}B(v_{\omega_{i}^{D}})\right)  \overset{S_{n-1}^{D}%
}{\rightarrow}B(v_{\Lambda_{n-1}^{D}})\otimes(B(v_{\Lambda_{n-1}^{D}})\cup
B(v_{\Lambda_{n}^{D}})).
\end{align*}

\item  Let $w$ be the reading of an admissible column $C$ of type $D$. If
$h(C)\prec n,$ denote by $t\ $the greatest unbarred letter such that
$t\notin\left|  C\right|  $. Write

- $l\frak{C}$ for the spin column of height $n$ obtained by adding to $lC$ the
barred letters $\overline{x}$ such that $x\notin\left|  C\right|  $.

- $r\frak{C}$ for the spin column of height $n$ obtained by adding to $rC$ the
unbarred letters $x$ such that $x\notin\left|  C\right|  $.

- $l_{t}\frak{C}$ for the spin column of height $n$ obtained by adding to $lC$
the letter $t$ and the barred letters $\overline{x}$ such that $x\notin\left|
C\right|  \cup\{t\}$.

- $r_{t}\frak{C}$ for the spin column of height $n$ obtained by adding to $rC$
the letter $\overline{t}$ and the unbarred letters $x$ such that
$x\notin\left|  C\right|  \cup\lbrack t\}$.
\end{enumerate}

\noindent Then we have
\[
\mathrm{(i)}:\left\{
\begin{tabular}
[c]{l}%
$S_{n}^{D}(w)=r\frak{C}\otimes l\frak{C}$ if $r\frak{C}\in B(v_{\Lambda
_{n}^{D}})$\\
$S_{n}^{D}(w)=r_{t}\frak{C}\otimes l_{t}\frak{C}$ otherwise
\end{tabular}
\right.  \text{ and }\mathrm{(ii)}:\left\{
\begin{tabular}
[c]{l}%
$S_{n-1}^{D}(w)=r\frak{C}\otimes l\frak{C}$ if $r\frak{C}\in B(v_{\Lambda
_{n-1}^{D}})$\\
$S_{n-1}^{D}(w)=r_{t}\frak{C}\otimes l_{t}\frak{C}$ otherwise
\end{tabular}
\right.
\]
(recall that $r\frak{C}\in B(v_{\Lambda_{n}^{D}})$ if and only if it contains
an even number of barred letters).
\end{lemma}

\begin{proof}
We only sketch the proof for $S_{n}^{D},$ the arguments are analogous for
$S_{n-1}^{D}.$

$1:$ The highest weight vertices of $B(v_{\Lambda_{n}^{D}})\otimes
(B(v_{\Lambda_{n}^{D}})\cup B(v_{\Lambda_{n-1}^{D}}))$ are the vertices
$v_{p}^{D}=\frak{C}_{n}\otimes\frak{C}_{p}$ with $p=0,...,n$. We have
$\mathrm{wt}(v_{p}^{D})=\omega_{p}^{D}$ for $p=1,...,n$ and $\mathrm{wt}%
(v_{0}^{D})=0.\;$Hence $S_{n}^{D}$ is the crystal isomorphism which sends
$B(v_{\omega_{p}^{D}})$ on $B(v_{p}^{D})$ for $p=1,...,n$ and $B(0)$ on
$v_{0}^{D}.$

$2:$When $w=v_{\omega_{p}^{D}}$, the equality $S_{n}^{D}(w)=r\frak{C}\otimes
l\frak{C}$ is true. Consider $w\in B(v_{\omega_{p}^{D}})$ and $i=1,...,n$ such
that $w^{\prime}=\widetilde{f}_{i}(w)\neq0$.\ Write $w=\mathrm{w}(C)$ and
$w^{\prime}=\mathrm{w}(C^{\prime})$ where $C$ and $C^{\prime}$ are two
admissible columns of height $p$.\ Let $t^{\prime}$ be the greatest unbarred
letter such that $t^{\prime}\notin\left|  C^{\prime}\right|  $. If the number
of barred letters of $C$ is equal to that of $C^{\prime},$ $r\frak{C}$ and
$r\frak{C}^{\prime}$ belongs together in $B(v_{\Lambda_{n}^{D}})$ or in
$B(v_{\Lambda_{n-1}^{D}}).$ In these cases we can prove that
\begin{equation}
S_{n}^{D}(w)=r\frak{C}\otimes l\frak{C\Longrightarrow} S_{n}^{D}(w^{\prime
})=r\frak{C}^{\prime}\otimes l\frak{C}^{\prime}\text{ and }S_{n}^{D}%
(w)=r_{t}\frak{C}\otimes l_{t}\frak{C\Longrightarrow} S_{n}^{D}(w^{\prime
})=r_{t^{\prime}}\frak{C}^{\prime}\otimes l_{t^{\prime}}\frak{C}^{\prime}
\label{f_SD}%
\end{equation}
as we have done for $S^{B}.$ Otherwise we have $i=n$ and $rC\cap E_{n}=(n-1\}$
or $rC\cap E_{n}=(n\}.$

\noindent Suppose $i=n$ and $n\in\left|  C\right|  $.\ Then $n-1$ is the
unique letter of $E_{n}=\{n-1,n,\overline{n},\overline{n-1}\}$ that occurs in
$C$.\ We have $t=n$ and $t^{\prime}=n-1$\ because $lC^{\prime}\cap
E_{n}=\overline{n}$.\ So $r\frak{C}\cap E_{n}=\{n,n-1\},$ $r_{t}\frak{C}\cap
E_{n}=\{\overline{n},n-1\},$ $l\frak{C}\cap E_{n}=\{\overline{n},n-1\}$ and
$l_{t}\frak{C}\cap E_{n}=\{n,n-1\}$.\ Similarly $r\frak{C}^{\prime}\cap
E_{n}=\{\overline{n},n-1\},$ $r_{t}\frak{C}^{\prime}\cap E_{n}=\{\overline
{n},\overline{n-1}\},$ $l\frak{C}^{\prime}\cap E_{n}=\{\overline{n}%
,\overline{n-1}\}$ and $l_{t}\frak{C}^{\prime}\cap E_{n}=\{\overline{n}%
,n-1\}$. Hence $\widetilde{f}_{i}(r\frak{C}\otimes l\frak{C)=}r_{t^{\prime}%
}\frak{C}^{\prime}\otimes l_{t^{\prime}}\frak{C}^{\prime}$ and $\widetilde
{f}_{i}(r_{t}\frak{C}\otimes l_{t}\frak{C)=}r\frak{C}^{\prime}\otimes
l\frak{C}^{\prime}$.\ We have
\begin{equation}
S_{n}^{D}(w)=r\frak{C}\otimes l\frak{C\Longrightarrow} S_{n}^{D}(w^{\prime
})=r_{t^{\prime}}\frak{C}^{\prime}\otimes l_{t^{\prime}}\frak{C}^{\prime
}\text{ and }S_{n}^{D}(w)=r_{t}\frak{C}\otimes l_{t}\frak{C\Longrightarrow}%
S_{n}^{D}(w^{\prime})=r\frak{C}^{\prime}\otimes l\frak{C}^{\prime}.
\label{f_SD1}%
\end{equation}
When $i=n$ and $n-1\in\left|  C\right|  $, we obtain (\ref{f_SD1}) by similar
arguments. Finally $\mathrm{(i)}$ follows from (\ref{f_SD}) and (\ref{f_SD1}).
\end{proof}

\begin{example}
Suppose $n=7$ and consider the admissible column $C$ of type $D$ such that
$\mathrm{w}(C)=67\overline{7}7\overline{6}.$ Then $\mathrm{w(}lC)=3457\bar{6}%
$, $\mathrm{w(}rC)=67\bar{5}\bar{4}\bar{3}.$ So $(t,\overline{t}%
)=(2,\overline{2})$ and, by identifying the spin columns with the set of
letters that they contain, we have $l\frak{C}=\{3457\bar{6}\bar{2}\bar{1}\},$
$r\frak{C}=\{1267\bar{5}\bar{4}\bar{3}\},$ $l_{t}\frak{C}=\{23457\bar{6}%
\bar{1}\}$, $r_{t}\frak{C}=\{167\bar{5}\bar{4}\bar{3}\bar{2}\}$. We have
$S_{n}^{D}(\mathrm{w}(C))=r_{t}\frak{C}\otimes l_{t}\frak{C}$ and $S_{n-1}%
^{D}(\mathrm{w}(C))=r\frak{C}\otimes l\frak{C}$ for $r\frak{C}\not \in
B(v_{\Lambda_{n}^{D}})$.
\end{example}

\begin{figure}
[ptb]
\begin{center}
\includegraphics[
height=6.4734cm,
width=10.1067cm
]%
{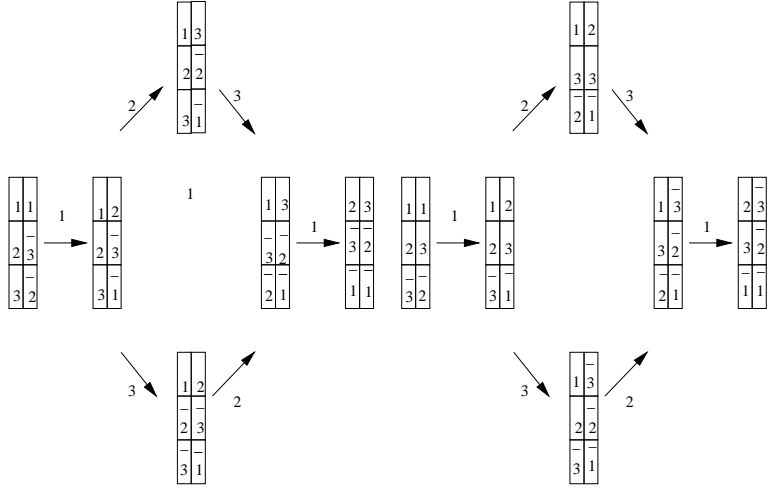}%
\caption{The connected components of $V(\Lambda_{3}^{D})^{\bigotimes2}$ and
$V(\Lambda_{2}^{D})^{\bigotimes2}$ isomorphic to $V(\omega_{1}^{D})$ for
$U_{q}(so_{6})$}%
\label{FIG2}%
\end{center}
\end{figure}

Although $C$ must be the empty column in Lemmas \ref{Lem_S_B} and
\ref{Lem_S_D}, we only use these Lemmas with $h(C)\geq1$ in the
sequel..\ Figure (\ref{FIG2}) below describe the connected components of
$V(\Lambda_{3}^{D})^{\otimes2}$ and $V(\Lambda_{2}^{D})^{\otimes2}$ isomorphic
to the vector representation $V(\Lambda_{1}^{D})$ of $U_{q}(so_{6})$ (see also
(\ref{vect_B})).%

Note that it is possible to describe explicitly the isomorphisms $(S^{B}%
)^{-1}$, $(S_{n}^{D})^{-1}$ and $(S_{n-1}^{D})^{-1}.\;$The reader interested
by this subject is referred to \cite{Lec2}.

\subsection{Plactic monoid for $\frak{G}_{n}\label{sub_sec_mon_spin}$}

Let $\lambda$ be a dominant weight such that $\lambda\notin\Omega_{+}$. If
$\lambda\in P_{+}^{B}$ then $\lambda$ has a unique decomposition
$\lambda=\Lambda_{n}^{B}+\lambda^{^{\prime}}$ with $\lambda^{\prime}\in
\Omega_{+}^{B}$. We set $v_{\lambda}^{B}=v_{\lambda^{\prime}}\otimes
v_{\Lambda_{n}^{B}}$. Then $v_{\lambda}^{B}$ is the highest weight vector of
$B(v_{\lambda}^{B}),$ a connected component of $\frak{G}_{n}^{B}$ isomorphic
to $B^{B}(\lambda)$. Denote by $Y(\lambda)$ the diagram obtained by adding a
K.N-diagram of height $n$ to $Y(\lambda^{\prime})$.

When$\ \lambda\in P_{+}^{D}$, $\lambda$ has a unique decomposition of type
$\lambda=\Lambda_{n}^{D}+\lambda^{\prime}$ with $\lambda^{\prime}\in\Omega
_{+}^{D}$ and $\overline{\omega}_{n}^{D}$ not appearing in $\lambda^{\prime}$
or $\lambda=\Lambda_{n-1}^{D}+\lambda^{\prime}$ with $\lambda^{\prime}%
\in\Omega_{+}^{D}$ and $\omega_{n}^{D}$ not appearing in $\lambda^{\prime}$.
According to this decomposition we set $v_{\lambda}^{D}=v_{\lambda^{\prime}%
}\otimes v_{\Lambda_{n}^{D}}$ or $v_{\lambda}=v_{\lambda^{\prime}}\otimes
v_{\Lambda_{n-1}^{D}}$. Then $v_{\lambda}^{D}$ is the highest weight vector of
$B(v_{\lambda}^{D}),$ a connected component of $\frak{G}_{n}^{D}$ isomorphic
to $B^{D}(\lambda)$. If $Y(\lambda^{\prime})=(Y^{\prime},\varepsilon)$ (see
\ref{Def_Y(lambda)}) with $\varepsilon\in\{-,0,+\},$ we set $Y(\lambda
)=(Y,\varepsilon)$ where $Y$ is the diagram obtained by adding a K.N diagram
of height $n$ to $Y^{\prime}$.

Given a tabloid $\tau$ and a spin column $\frak{C}$, the spin tabloid
$[\frak{C},T]$ is obtained by adding $\frak{C}$ in front of $\tau$. The
reading of the spin tabloid $[\frak{C},\tau]$ is \textrm{w(}$[\frak{C}%
,\tau])=\mathrm{w(}\tau)\otimes\frak{C=}\mathrm{w(}\tau)\frak{C}$. Note that
the vertices of $B(v_{\lambda})$ are readings of spin tabloids.

\begin{definition}
\label{def_tab_spin} \ \ \ \ 

\begin{itemize}
\item  Let $\lambda\in P_{+}^{B}$ such that $\lambda\notin\Omega_{+}^{B}.$ A
spin tabloid is a spin tableau of type $B$ and shape $Y(\lambda)$ if its
reading is a vertex of $B(v_{\lambda}^{B})$.

\item  Let $\lambda\in P_{+}^{D}$ such that $\lambda\notin\Omega_{+}^{D}.$ A
spin tabloid is a spin tableau of type $D$ and shape $Y(\lambda)$ if its
reading is a vertex of $B(v_{\lambda}^{D})$.
\end{itemize}
\end{definition}

It follows from this definition that for $\frak{T}_{1}$ and $\frak{T}_{2}$ two
spin tableaux $\frak{T}_{1}\sim\frak{T}_{2}\Longleftrightarrow\frak{T}%
_{1}=\frak{T}_{2}$. It is possible to extend Definition \ref{Def_b_conf} to a
spin tableau $[\frak{C,}C]$ of type $D$ with $C$ an admissible column of type
$D$. We will say that $[\frak{C,}C]$ contains an $a$-configuration even or odd
when this configuration appears in the tableau of two columns $C_{\frak{C}}C$
where $C_{\frak{C}}$ is the admissible column of type $D$ and height $n$
containing the letters of $\frak{C}.$ Kashiwara and Nakashima have obtained in
\cite{KN} a combinatorial description of the orthogonal spin tableaux
equivalent to the following:

\begin{theorem}
\label{TH_KNS} \ \ \ \ \ 

\begin{itemize}
\item $\frak{T=}[\frak{C,}T]$ is a spin tableau of type $B$ if and only if $T$
is a tableau of type $B$ and the rows of $[\frak{C,}lC_{1}]$ weakly increase
from left to right.

\item $\frak{T=}[\frak{C,}T]$ is a spin tableau of type $D$ if and only if $T$
is a tableau of type $D$, the rows of $[\frak{C,}lC_{1}]$ weakly increase from
left to right and $[\frak{C,}lC_{1}]$ does not contain an $a$-configuration
(even or odd) with $q(a)=n-a$.
\end{itemize}
\end{theorem}

\noindent It follows from the definition above that for any spin tableau
$[\frak{C,}T]$ of type $D$%

\begin{align*}
\frak{C}  &  \in B(\Lambda_{n}^{D})\text{ implies that the shape of }T\text{
is }(Y,\varepsilon)\text{ with }\varepsilon\neq-,\\
\frak{C}  &  \in B(\Lambda_{n-1}^{D})\text{ implies that the shape of }T\text{
is }(Y,\varepsilon)\text{ with }\varepsilon\neq+\text{.}%
\end{align*}
A generalized tableau is an orthogonal tableau or a spin orthogonal tableau.
Similarly to subsection \ref{subsec_monoids}, the quotient sets $\frak{G}%
_{n}/\overset{B}{\sim}$ and $\frak{G}_{n}/\overset{D}{\sim}$ can be
respectively identified with the sets of generalized tableaux of type $B$ and
$D$. For $x$ a letter of $\mathcal{B}_{n}$ or $\mathcal{D}_{n}$ and $\frak{C}$
a spin column of height $n$ whose greatest letter is $z$, we write
$x\vartriangle\frak{C}$ when $x\nleq z.$

\begin{definition}
\label{def_mono_spinB}The monoid $\frak{Pl(}B_{n}\frak{)}$ is the quotient set
of $\frak{B}_{n}^{\ast}$ by the relations:

\begin{itemize}
\item $R_{i}^{B},$ $i=1,...,5$ defining $Pl(B_{n})$,

\item $R_{6}^{B}$: for $x\in\mathcal{B}_{n}$ and $\frak{C}$ a spin column such
that $x\vartriangle\frak{C;}$ $\frak{C}x\equiv\frak{C}^{\prime}$ where
$\frak{C}^{\prime}$ is the spin column such that $\mathrm{wt(}\frak{C}%
^{\prime})=\mathrm{wt(}\frak{C})+\mathrm{wt(}x),$

\item $R_{7}^{B}$: for $x\in\mathcal{B}_{n}$ and $\frak{C}$ a spin column such
that $x\not \vartriangle\frak{C;}$ $\frak{C}x\equiv x^{\prime}\frak{C}%
^{\prime}$ where
\[
\left\{
\begin{tabular}
[c]{l}%
$x^{\prime}=\min\{t\in\frak{C};$ $t\succeq x\}$ if $x\succeq0$\\
$x^{\prime}=\min\{t\in\frak{C};$ $t\succeq x\}\cup\{0\}$ if $x\preceq n$%
\end{tabular}
\right.
\]
and $\frak{C}^{\prime}$ is the spin column such that $\mathrm{wt}%
(\frak{C}^{\prime})=\mathrm{wt(}\frak{C})+\mathrm{wt(}x)-\mathrm{wt(}%
x^{\prime})$,

\item $R_{8}^{B}$:for $C$ an admissible column of type $B,$ $S^{B}%
(\mathrm{w}(C))\equiv\mathrm{w}(C)$.
\end{itemize}
\end{definition}

Lemma \ref{lem_phi_tens} implies that the highest weight vertex of the
connected component containing a word $\frak{C}x$ with $x\in\mathcal{B}_{n}$
and $\frak{C}$ a spin column may be written $\frak{C}_{n}x_{0}$ where
$x_{0}\in\{0,1\}$. So $\frak{C}x\in B(v_{\Lambda_{n}^{B}}\otimes0)$ or
$\frak{C}x\in B(v_{\Lambda_{n}^{B}}\otimes1).\;$The following lemma gives the
interpretation of relations $R_{6}^{B}$ and $R_{7}^{B}$ in terms of crystal isomorphisms.

\begin{lemma}
\label{Lem_iso_spin} \ \ \ \ \ \ \ \ \ \ 

\begin{enumerate}
\item  The vertices of $B(v_{\Lambda_{n}^{B}}\otimes0)$ are the words of the
form $\frak{C}x$ where $\frak{C}$ is a spin column and $x\in\mathcal{B}_{n}$
such that $x\vartriangle\frak{C.}$

\item  The vertices of $B(v_{\Lambda_{n}^{B}}\otimes1)$ are the words of the
form $\frak{C}x$ where $\frak{C}$ is a spin column and $x\in\mathcal{B}_{n}$
such that $x\not \vartriangle\frak{C.}$

\item  Denote by $\Psi$ and $\Psi^{\prime}$ the crystal isomorphisms:
\[
\text{%
\begin{tabular}
[c]{l}%
$\Psi:B(v_{\Lambda_{n}^{B}}\otimes0)\rightarrow B(v_{\Lambda_{n}^{B}})$\\
$\Psi^{\prime}:B(v_{\Lambda_{n}^{B}}\otimes1)\rightarrow B(1\otimes
v_{\Lambda_{n}^{B}})$%
\end{tabular}
.}%
\]
Then if the word $\frak{C}x$ occur in the left hand side a relation $R_{6}%
^{B}$ (resp. of $R_{7}^{B}),$ $\Psi(\frak{C}x)$ (resp. $\Psi^{\prime}%
(\frak{C}x))$ is the word occurring in the right hand side of this relation.
\end{enumerate}
\end{lemma}

\begin{proof}
$1$ Consider a word $\frak{C}x$ such that $x\vartriangle\frak{C}$ and
$\widetilde{f}_{i}(\frak{C}x)\neq0$.\ Let $y$ be the greatest letter of
$\frak{C}$.\ Set $\widetilde{f}_{i}(\frak{C}x)=\frak{U}t$ where $\frak{U}$ is
a spin column and $t$ a letter of $\mathcal{B}_{n}$. We are going to show that
$t\vartriangle\frak{U}$.\ If $y$ is the greatest letter of $\frak{U}$ then
$t\succeq x\succ y,$ hence $t\vartriangle\frak{U}$.\ Otherwise $\widetilde
{f}_{i}(\frak{C}x)=\widetilde{f}_{i}(\frak{C)}x$ thus $\varepsilon_{i}(x)=0$
by (\ref{TENS1})$.\;$When $i\neq n,$ we must have $y=\overline{i+1}$, $x\succ
y$ and $x\notin\{\overline{i},i+1\}$ because $\varepsilon_{i}(x)=0$.\ Hence
$x\succ\overline{i}$ and $x=t\vartriangle\frak{U}$ for $\overline{i}$ is the
greatest letter of $\frak{U}$.\ When $i=n,$ $y=n$ and $x\succ\overline{n}$
because $\varepsilon_{n}(x)=0$.\ We obtain similarly $t\vartriangle\frak{U}%
$..\ Hence the set of words $\frak{C}x$ such that $x\vartriangle\frak{C}$ is
closed under the action of the $\widetilde{f}_{i}.$ By similar arguments we
can prove that this set is also closed under the action of the $\widetilde
{e}_{i}.$ Moreover $v_{\Lambda_{n}^{B}}\otimes0$ is the unique highest weight
vertex among these words $\frak{C}x$. Hence $B(v_{\Lambda_{n}^{B}}\otimes0)$
contains exactly the words of the form $\frak{C}x$ such that $x\vartriangle
\frak{C}$.\ 

$2$ Follows immediately from $1.$

$3$ If $x\vartriangle\frak{C}$, $\Psi(\frak{C}x)$ is the unique spin column of
weight $\mathrm{wt(}\frak{C}x)$, that is $\Psi(\frak{C}x)=\frak{C}^{\prime}$
with the notation of $R_{6}^{B}$. When $x\not \vartriangle\frak{C,}$ we
consider the following cases:

$\mathrm{(i)}$: $x\in\frak{C}$.\ Set $\Psi(\frak{C}x)=y\frak{D}$. Then we
deduce from the equality $\mathrm{wt(}y\frak{D})=\mathrm{wt(}\frak{C}x)$ that
$y=x$ and $\frak{D=C}$. Indeed $x\frak{C}$ is the unique vertex of
$B(1)\otimes B(v_{\Lambda_{n}^{B}})$ of weight $\mathrm{wt}(\frak{C}x).$ Hence
$y=x=t$ and $\frak{D=C}^{\prime}$ with the notation of $R_{6}^{B}$.

$\mathrm{(ii)}$: $x\notin\frak{C}$ .$\;$When $x\succ0,$ set $x=\overline{p}$
and $\overline{k}=\min\{t\in\frak{C};$ $t\succeq x\}$. Then
$\{p,p-1,...,k+1\}\subset\frak{C}$. By using the formulas (\ref{TENS1}) and
(\ref{TENS2}) we obtain
\[
\widetilde{f}_{k}\cdot\cdot\cdot\widetilde{f}_{p-2}\widetilde{f}%
_{p-1}(\frak{C}\overline{p})=\frak{C}\overline{k}%
\]
So, by $\mathrm{(i)}$, $\frak{C}\overline{k}\sim\overline{k}\frak{C}$ which
implies
\[
\frak{C}\overline{p}\sim\widetilde{e}_{p-1}\cdot\cdot\cdot\widetilde{e}%
_{k}(\overline{k}\frak{C})=\overline{k}\widetilde{e}_{p-1}\cdot\cdot
\cdot\widetilde{e}_{k}(\frak{C})=\overline{k}\frak{C}^{\prime}%
\]
with the notation of $R_{7}^{B}$. It means that $\Psi(\frak{C}x)=\overline
{k}\frak{C}^{\prime}$. When $x=0,$ we have $\widetilde{f}_{x^{\prime}-1}%
\cdot\cdot\cdot\widetilde{f}_{1}\widetilde{f}_{n}(\frak{C}0)=\frak{C}%
\overline{k}$.because $\{n,n-1,...,k+1\}\subset\frak{C}$ and we terminate as
above. When $x=p\prec0$ and $\min\{t\in\frak{C};$ $t\succeq x\}\cup
\{0\}=k\prec0$, we have $\{\overline{p},\overline{p+1},...,\overline
{k-1}\}\subset\frak{C}.$ So $\widetilde{f}_{k-1}\cdot\cdot\cdot\widetilde
{f}_{p+1}\widetilde{f}_{p}(\frak{C}p)=\frak{C}k$ and the proof is similar. If
$\min\{t\in\frak{C};$ $t\succeq p\}\cup\{0\}=0$, $\{\overline{p}%
,\overline{p+1},...,\overline{n}\}\subset\frak{C}$. Then $\widetilde{f}%
_{n}\cdot\cdot\cdot\widetilde{f}_{p+1}\widetilde{f}_{p}(\frak{C}%
p)=\frak{C}0\sim\overline{n}\frak{C}^{%
{{}^\circ}%
}$ with $\frak{C}^{%
{{}^\circ}%
}=\frak{C}-\{\overline{n}\}+\{n\}$ by the case $x=0$.\ So formulas
(\ref{TENS1}) and (\ref{TENS2}) imply that $\frak{C}x\sim\widetilde{e}%
_{p}\cdot\cdot\cdot\widetilde{e}_{n}(\overline{n}\frak{C}^{%
{{}^\circ}%
})=\widetilde{e}_{n}(\overline{n})\widetilde{e}_{p}\cdot\cdot\cdot
\widetilde{e}_{n-1}(\frak{C}^{%
{{}^\circ}%
})=0\frak{C}^{\prime}$ with the notation of $R_{7}^{B}$. It means that
$\Psi(\frak{C}x)=0\frak{C}^{\prime}$.
\end{proof}

\begin{definition}
\label{def_mono_spinD}The monoid $\frak{Pl(}D_{n}\frak{)}$ is the quotient set
of $\frak{D}_{n}^{\ast}$ by the relations:

\begin{itemize}
\item $R_{i}^{D},$ $i=1,...,5$ defining $Pl(D_{n})$,

\item $R_{6}^{D}$: for $x\in\mathcal{D}_{n}$ and $\frak{C}$ a spin column such
that $x\vartriangle\frak{C;}$ $\frak{C}x\equiv\frak{C}^{\prime}$ where
$\frak{C}^{\prime}$ is the spin column such that $\mathrm{wt(}\frak{C}%
^{\prime})=\mathrm{wt(}\frak{C})+\mathrm{wt(}x)$,

\item $R_{7}^{D}$: for $x\in\mathcal{D}_{n}$ and $\frak{C}$ a spin column such
that $x\not \vartriangle\frak{C;}$ $\frak{C}x\equiv x^{\prime}\frak{C}%
^{\prime}$ where $x^{\prime}=\min\{t\in\frak{C};$ $t\succeq x\}$ and
$\frak{C}^{\prime}$ is the spin column such that $\mathrm{wt(}\frak{C}%
^{\prime})=\mathrm{wt(}\frak{C})+\mathrm{wt(}x)-\mathrm{wt(}x^{\prime})$,

\item $R_{8}^{D}$: for $C$ an admissible column of type $D,$ $S_{n}%
^{D}(\mathrm{w}(C))\equiv\mathrm{w}(C)$ and $S_{n-1}^{D}(\mathrm{w}%
(C))\equiv\mathrm{w}(C).$
\end{itemize}
\end{definition}

We can prove by using similar arguments to those of Lemma \ref{Lem_iso_spin}
that the relations $R_{6}^{D}$ and $R_{7}^{D}$ read from left to right
describe respectively the crystal isomorphisms
\begin{equation}%
\begin{tabular}
[c]{l}%
$\left\{
\begin{tabular}
[c]{l}%
$B(v_{\Lambda_{n}^{D}}\otimes\overline{n})\rightarrow B(v_{\Lambda_{n-1}^{D}%
})$\\
$B(v_{\Lambda_{n-1}^{D}}\otimes n)\rightarrow B(v_{\Lambda_{n}^{D}})$%
\end{tabular}
\right.  $\\
and\\
$\left\{
\begin{tabular}
[c]{l}%
$B(v_{\Lambda_{n}^{D}}\otimes1)\rightarrow B(1\otimes v_{\Lambda_{n}^{D}})$\\
$B(v_{\Lambda_{n-1}^{D}}\otimes1)\rightarrow B(1\otimes v_{\Lambda_{n-1}^{D}%
})$%
\end{tabular}
\right.  .$%
\end{tabular}
\label{interpret_R7D}%
\end{equation}

\begin{lemma}
\label{lem_copat_cry_op}Let $w_{1}$ and $w_{2}$ be\ two vertices of
$\frak{G}_{n}$ such that $w_{1}\equiv w_{2}$. Then for $i=1,...,n$:
\begin{align*}
\widetilde{e}_{i}(w_{1})  &  \equiv\widetilde{e}_{i}(w_{2})\text{ and
}\varepsilon_{i}(w_{1})=\varepsilon_{i}(w_{2}),\\
\widetilde{f}_{i}(w_{1})  &  \equiv\widetilde{f}_{i}(w_{2})\text{ and }%
\varphi_{i}(w_{1})=\varphi_{i}(w_{2}).
\end{align*}
\end{lemma}

\begin{proof}
By induction we can suppose that $w_{2}$ is obtained from $w_{1}$ by applying
only one plactic relation. In this case we write $w_{1}=u\widehat{w}_{1}v$ and
$w_{2}=u\widehat{w}_{2}v$ where $u,v,\widehat{w}_{1},\widehat{w}_{2}$ are
factors of $w_{1}\ $and $w_{2}$ such that $\widehat{w}_{1}\equiv\widehat
{w}_{2}$ by one of the relations $R_{i}$. Formulas (\ref{TENS1}) and
(\ref{TENS2}) imply that it is enough to prove the lemma for $\widehat{w}_{1}
$ and $\widehat{w}_{2}$. This last point is immediate because we have seen
that each plactic relation may be interpreted in terms of a crystal isomorphism.
\end{proof}

So we obtain $w_{1}\equiv w_{2}\Longrightarrow w_{1}\sim w_{2}.$ To establish
the implication $w_{1}\sim w_{2}\Longrightarrow w_{1}\equiv w_{2}$, it
suffices, as in subsection \ref{subsec_monoids} to prove that two highest
weight vertices of $\frak{G}_{n}^{B}$ (resp.\ $\frak{G}_{n}^{D}$) with the
same weight are congruent in $\frak{Pl(}B_{n}\frak{)}$ (resp.\ $\frak{Pl(}%
D_{n}\frak{)}$). Given a vertex $w\in\frak{G}_{n},$ we know by Theorems
\ref{TH_KNS} and \ref{TH_KN} that there exists a unique generalized tableau
$\frak{P}(w)$ such that
\[
\mathrm{w(}\frak{P}(w))\sim w.
\]

\begin{lemma}
Let $w$ be a highest weight vertex of $\frak{G}_{n}$. Then $\mathrm{w(}%
\frak{P}(w))\equiv w$.
\end{lemma}

\begin{proof}
By using relations $R_{6}$ and $R_{7}$, $w$ is congruent to a word $u\frak{U}
$ such that $u\in G_{n}$ and $\frak{U\in G}_{n}$.\ Relation $R_{8}$ implies
that any word consisting in an even number of spin columns is congruent to a
vertex of $G_{n}$.\ If $\frak{U}$ contains an even number of spin columns,
there exists $v\in G_{n}$ such that $w\equiv v$.\ We have $\frak{P}(w)=P(v)$
because $w\equiv v\Longrightarrow w\sim v.$ Thus $\mathrm{w(}\frak{P}%
(w))=\mathrm{w(}P(v))\equiv v\equiv w$ and the lemma is proved. If $w$
contains an odd number of spin columns, there exists a vertex $v\in G_{n}$ and
a spin column $\frak{C}$ such that $w\equiv v\frak{C}$.\ Set $P(v)=T$.\ Then
$w\equiv\mathrm{w}(T)\frak{C}$. Write $T=C\widehat{T}$ where $C$ is the first
column of $T$ and $\widehat{T}$ the tableau obtained by erasing $C$ in $T$. By
Lemma \ref{lem_phi_tens}, $\mathrm{w}(T)$ is a highest weight vertex because
$w$ is a highest weight vertex of $\frak{G}_{n}$.\ In particular,
$\mathrm{w}(C)$ is a highest weight vertex. Set $p=h(C)$.

Suppose first $w\in\frak{G}_{n}^{B}$.\ We have $S^{B}(\mathrm{w}%
(C))=\frak{C}_{n}\frak{C}_{p}$ (see Lemma \ref{Lem_S_B}). So $w\equiv
\mathrm{w}(\widehat{T})\frak{C}_{n}\frak{C}_{p}\frak{C}$.\ By Lemma
\ref{lem_phi_tens} we must have $\varepsilon_{i}(\frak{C})=0$ for
$i=p+1,...,n$. This implies that the letters of $\{\overline{p+1}%
,...,\overline{n}\}$ do not appear in $\frak{C}$. Indeed $\overline{n}%
\notin\frak{C}$ otherwise $\varepsilon_{n}(\frak{C})\neq0$ and if
$\overline{q}\succ\overline{n}$ is the lowest barred letter of $\{\overline
{p+1},...,\overline{n}\}$ appearing in $\frak{C} $ we obtain $\varepsilon
_{q}(\frak{C})=1\neq0$ because $q+1\in\frak{C}$. So $\frak{C}$ contains the
letters of $\{p+1,...,n\}$. Let $\{x_{1}\prec\cdot\cdot\cdot\prec x_{s}\}$ be
the set of unbarred letters $\preceq p$ that occur in $\frak{C}$.\ By Lemma
\ref{Lem_S_B}, we have
\[
S^{B}(x_{1}\cdot\cdot\cdot x_{s}\underset{n-p\text{ times}}{\underbrace
{0\cdot\cdot\cdot0}})=\frak{C}_{p}\frak{C.}%
\]
Hence
\[
w\equiv\mathrm{w(}\widehat{T})\frak{C}_{n}(x_{1}\cdot\cdot\cdot x_{s}%
\underset{n-p\text{ times}}{\underbrace{0\cdot\cdot\cdot0}})
\]
and by applying relations $R_{6}^{B}$ and $R_{7}^{B}$ we have $w\equiv
\mathrm{w(}\widehat{T})(x_{1}\cdot\cdot\cdot x_{s})\frak{C}_{n}$.\ Write
$T^{\prime}=x_{s}\rightarrow(\rightarrow\cdot\cdot\cdot x_{1}\rightarrow
\widehat{T})$. Then $[\frak{C}_{n},T^{\prime}]$ is a spin orthogonal tableau
and $\mathrm{w(}T^{\prime})\frak{C}_{n}\equiv w$. So $T^{\prime}=\frak{P}(w)$
and the lemma is true.

Suppose now $w\in\frak{G}_{n}^{D}$.\ If the shape of $\widehat{T}$ is
$(Y,\varepsilon)$ with $\varepsilon\not =-,$ we consider $S_{n}^{D}%
(\mathrm{w}(C))=\frak{C}_{n}\frak{C}_{p}$. Then $[\frak{C}_{n},\widehat{T}]$
is a spin tableau and the proof is similar to that of the type $B$ case. If
the shape of $\widehat{T}$ is $(Y,\varepsilon)$ with $\varepsilon=-,$ it
suffices to consider $S_{n-1}^{D}(\mathrm{w}(C))=\frak{C}_{n-1}\frak{C}_{n-1}$
where instead of $S_{n}^{D}(\mathrm{w}(C))$.
\end{proof}

Now if $w_{1}$ and $w_{2}$ are two highest weight vertices of $\frak{G}_{n}$
with the same weight $\lambda$, we have $\frak{P}(w_{1})=\frak{P}(w_{2})$
because there is only one orthogonal tableau of highest weight $\lambda
$.\ Then the lemma above implies that $w_{1}\equiv w_{2}$. We can state the

\begin{theorem}
Let $w_{1}$ and $w_{2}$ two vertices of $\frak{G}_{n}\frak{.}$ Then $w_{1}\sim
w_{2}$ if and only if $w_{1}\equiv w_{2}$.
\end{theorem}

For any vertex $w\in\frak{G}_{n}$, it is possible to obtain $\frak{P}%
(w)\frak{\ }$by using an insertion algorithm analogous to that describe in
Section \ref{sec_in_vect}.\ Considering the sequence of shape of the
intermediate generalized tableaux appearing during the computation of
$\frak{P}(w)$, we obtain a $\frak{Q}$-symbol $\frak{Q}(w)$. Then for $w_{1}$
and $w_{2}$ two vertices of $\frak{G}_{n}$ we have:
\[
w_{1}\longleftrightarrow w_{2}\Longleftrightarrow\frak{Q}(w_{1})=\frak{Q}%
(w_{1})
\]
where $w_{1}\longleftrightarrow w_{2}$ means that $w_{1}$ and $w_{2}$ occur in
the same connected component of $\frak{G}_{n}$. The reader interested by this
subject is referred to \cite{Lec2}.

\end{document}